\documentclass[12pt]{article}
\usepackage{amsmath}
\usepackage{amssymb}
\numberwithin{equation}{section}
\usepackage{epic}

\usepackage{graphicx}

\title{Universality of the Pearcey process}

\author{ Mark Adler\thanks
{2000 {\em Mathematics Subject Classification}. Primary:
60J60, 60J65, 60G55; secondary: 35Q53, 35Q58. {\em Key
words and Phrases}: Non-intersecting Brownian motions,
Pearcey distribution, matrix models, random Hermitian
ensembles, multi-component KP equation, Virasoro
constraints.
\newline
  Department of Mathematics, Brandeis University,
Waltham, MA 02454, USA. E-mail: adler@brandeis.edu. The
support of a National Science Foundation grant \#
DMS-04-06287 is gratefully acknowledged.}~~~~~~
Nicolas Orantin\thanks{Universit\'e de Louvain, 1348 Louvain-la-Neuve, Belgium. E-mail: orantin.nicolas@uclouvain.be.
The support of the ENRAGE European network MRTN-CT-2004-005616,
the ENIGMA European network MRT-CT-2004-5652, the French and Japanese governments through PAI SAKURA,
the European Science Foundation through the MISGAM program and the ANR project G\'{e}om\'{e}trie Int\'{e}grabilit\'{e}
en Physique Math\'{e}matique ANR-BLAN-0029-01 is gratefully acknowledged.}
  ~~~~~~ Pierre
van Moerbeke\thanks{ Department of Mathematics,
Universit\'e de Louvain, 1348 Louvain-la-Neuve, Belgium
and Brandeis University, Waltham, MA 02454, USA. E-mail:
vanmoerbeke@math.ucl.ac.be and @brandeis.edu. The
support of a National Science Foundation grant \#
DMS-04-06287, a European Science Foundation grant
(MISGAM), a Marie Curie Grant (ENIGMA), Nato, FNRS and
Francqui Foundation grants is gratefully acknowledged.}}

\date{}

\newcommand{\MAT}[1]{\left(\begin{array}{*#1c}}
\newcommand{\mat}{\end{array}\right)}
\newcommand{\qed}{\leavevmode\unskip\nobreak\penalty200\hskip2pt\null
\nobreak\hfill\rule{1.1ex}{1.1ex}
\medbreak }

\newcommand{\rg}{\rightarrow}

\newcommand{\CR}{{\cal C}}

\newcommand{\HR}{{\cal H}}

\newcommand{\PR}{{\cal P}}

\newcommand{\BP}{{\mathbb P}}

\newcommand{\BX}{{\mathbb X}}

\newcommand{\iy}{\infty}
\newcommand{\pl}{\partial}
\def\a{\alpha}
\def\b{\beta}
\def\e{\eta}
\def\t{\tau}

\newcommand{\no}{\nonumber}

\newcommand{\beq}{\begin{equation}}
\newcommand{\eeq}{\end{equation}}
\newcommand{\bea}{\begin{eqnarray}}
\newcommand{\eea}{\end{eqnarray}}

\newenvironment{proof}{\medskip\noindent{\it Proof:\/} }{\qed}
\newenvironment
         {remark}{\medskip\noindent\underline{\it Remark:\/} }{\medbreak}

\newcommand{\om}{\omega}

\newcommand{\la}{\langle}
\newcommand{\ra}{\rangle}
\newcommand{\ga}{\gamma}

\newcommand{\dt}{\delta}
\newcommand{\Dt}{\Delta}
  \newcommand{\vr}{\varepsilon}

\newcommand{\BR}{{\mathbb R}}
\newcommand{\lb}{\lambda}

\newcommand{\diag}{\operatorname{diag}}

\def\be#1\ee{\begin{equation}#1\end{equation}}
\def\bea#1\eea{\begin{eqnarray}#1\end{eqnarray}}
\def\bean#1\eean{\begin{eqnarray*}#1\end{eqnarray*}}

\newcommand{\Tr}{\operatorname{\rm Tr}}

\newtheorem{definition}{Definition}[section]
\newtheorem{theorem}[definition]{Theorem}
\newtheorem{lemma}[definition]{Lemma}
\newtheorem{corollary}[definition]{Corollary}
\newtheorem{proposition}[definition]{Proposition}

\def\br{\begin{remark}\rm\small}
\def\er{\end{remark}}
\def\bt{\begin{theorem}}
\def\et{\end{theorem}}
\def\bd{\begin{definition}}
\def\ed{\end{definition}}
\def\bp{\begin{proposition}}
\def\ep{\end{proposition}}
\def\bl{\begin{lemma}}
\def\el{\end{lemma}}
\def\bc{\begin{corollary}}
\def\ec{\end{corollary}}
\def\beaq{\begin{eqnarray}}
\def\eeaq{\end{eqnarray}}

\newcommand{\eq}[1]{Eq.~(\ref{#1})}

\catcode `!=11

\newdimen\squaresize
\newdimen\thickness
\newdimen\Thickness
\newdimen\ll! \newdimen \uu! \newdimen\dd! \newdimen \rr! \newdimen
\temp!

\def\sq!#1#2#3#4#5{%
\ll!=#1 \uu!=#2 \dd!=#3 \rr!=#4
\setbox0=\hbox{%
  \temp!=\squaresize\advance\temp! by .5\uu!
  \rlap{\kern -.5\ll!
  \vbox{\hrule height \temp! width#1 depth .5\dd!}}%
  \temp!=\squaresize\advance\temp! by -.5\uu!
  \rlap{\raise\temp!
  \vbox{\hrule height #2 width \squaresize}}%
  \rlap{\raise -.5\dd!
  \vbox{\hrule height #3 width \squaresize}}%
  \temp!=\squaresize\advance\temp! by .5\uu!
  \rlap{\kern \squaresize \kern-.5\rr!
  \vbox{\hrule height \temp! width#4 depth .5\dd!}}%
  \rlap{\kern .5\squaresize\raise .5\squaresize
  \vbox to 0pt{\vss\hbox to 0pt{\hss $#5$\hss}\vss}}%
}
  \ht0=0pt \dp0=0pt \box0
}

\def\vsq!#1#2#3#4#5\endvsq!{\vbox to \squaresize{\hrule
width\squaresize height 0pt%
\vss\sq!{#1}{#2}{#3}{#4}{#5}}}

\newdimen \LL! \newdimen \UU! \newdimen \DD! \newdimen \RR!

\def\vvsq!{\futurelet\next\vvvsq!}
\def\vvvsq!{\relax
   \ifx     \next l\LL!=\Thickness \let\continue=\skipnexttoken!
   \else\ifx\next u\UU!=\Thickness \let\continue=\skipnexttoken!
   \else\ifx\next d\DD!=\Thickness \let\continue=\skipnexttoken!
   \else\ifx\next r\RR!=\Thickness \let\continue=\skipnexttoken!
   \else\def\continue{\vsq!\LL!\UU!\DD!\RR!}%
   \fi\fi\fi\fi
   \continue}

\def\skipnexttoken!#1{\vvsq!}

\def\place#1#2#3{\vbox to 0pt{\vss
\rlap{\kern#1\squaresize
   \raise#2\squaresize\hbox{$#3$}}
\vss}}

\catcode `!=12

\begin{document}
\maketitle

\tableofcontents

\begin{abstract}
Consider non-intersecting Brownian motions on the line leaving from the origin   and forced to two arbitrary points. Letting the number of Brownian particles tend to infinity, and upon rescaling, there is a point of bifurcation, where the support of the density of particles goes from one interval to two intervals. In this paper, we show that at that very point of bifurcation a cusp appears, near which the Brownian paths fluctuate like the Pearcey process. This is a universality result within this class of problems. Tracy and Widom obtained such a result in the symmetric case, when the two target points are symmetric with regard to the origin. This asymmetry enabled us to improve considerably a result concerning the non-linear partial differential equations governing the transition probabilities for the Pearcey process, obtained by Adler and van Moerbeke.

\end{abstract}
\section{Introduction}

Consider the probability that $n$ non-intersecting
(Dyson) Brownian motions
$$ x_1(t)<\ldots <x_n(t) $$ in $\BR$ belong to a set $E\in \BR$, with all particles
leaving from the origin at time $t=0$ and all forced to
end up at $b_1<b_2<\ldots<b_p$ at time $t=1$:

\bean   
\BP^{(b_1,\ldots,b_p)}_{n} \left(\begin{tabular}{c|c}
& all $x_j(0) =0$\\
all $x_j(t)\in E$ for $1\leq j\leq n$ &
 $n_1$ paths end up at $b_1$ at $t=1$\\
  & \vdots\\
  &
$n_p$ paths end up at $b_p$ at $t=1$
\end{tabular}\right)
\eean
with
$\sum_{i=1}^p n_i=n$ and with (local) transition
probability
 \be
 p(t;x,y):=\frac{1}{\sqrt{\pi t}} e^{-\frac{(x-y)^2}{t}}
 .\label{tr-pr}
  \ee
 A formula by Karlin-McGregor enables one to express this probability as an integral of a product of two determinants involving the transition probability (\ref{tr-pr}) above. This further leads to a expression as (i) a GUE-matrix integral with an external potential, (ii) a determinant of a block moment matrix, with $p$ blocks and (iii) a Fredholm determinant of a kernel. Finally it is also the solution of a PDE in the end-points of the interval $E$ and the target points $b_1,\ldots,b_p$.

Throughout this paper, we shall be dealing with the case of two target points $p=2$. In this paper, we show that, when $n\rg \iy$ and when one looks through a microscope near a certain point of bifurcation, the non-intersecting Brownian motions tend to a new process, the {\em Pearcey process}, whatever be the location of the target points and whatever be the proportion of particles forced to those points. Tracy and Widom \cite{TW-Pearcey} showed this result in the symmetric case; namely when the target points are symmetric with respect to the origin and half of the particles go to either target point. Br\'ezin and Hikami \cite{Brezin2,Brezin3,Brezin4,Brezin5} first considered this kernel and Bleher-Kuijlaars \cite{BleKui2} obtained strong asymptotics using Riemann-Hilbert techniques.

 The {\em Pearcey process} ${\cal P}(t)$ describes a cloud of
Brownian particles, evolving in time according to a
(matrix) Fredholm determinant,
$$ \BP^{\PR}\left(\mbox{all}~~{\cal P}(t_j)\in
E_j^c,1\leq j\leq m\right)=\det\left(I-\left(
\raisebox{1mm}{$\chi$}{}_{_{E_i}}
K^{\PR}_{t_it_j}\raisebox{1mm}{$\chi$}{}_{_{E_j}}
\right)_{1\leq i,j\leq m}\right) \label{2-time}$$
of the Pearcey kernel \bea{
 K^{\PR}_{s,t}(x,y)
 }
 &=&
- \frac{1
       }{4\pi^2 }\int_X dV
\int^{i\iy}_{-i\iy}dU 
e^{-\frac{U^4}{4}+\frac{tU^2}{2}-Uy}
 e^{\frac{V^4}{4}-\frac{sV^2}{2}+Vx}
\frac{1}{U-V} \no\\ \no\\
&&
%
-\frac{{\mathbb I}(s<t)}{\sqrt{2\pi (t-s)}} e^{-{\frac { \left( { x}-{
y} \right) ^{2}}{2 ({ t}-
 {
s})}}}
 \no\eea\be\label{Pkernel} \ee
The contour $X$ is given by the ingoing rays from $\pm
\iy e^{i\pi/4}$ to $0$ and the outgoing rays from $0$ to
 $\pm
\iy e^{-i\pi/4}$, i.e., $X$ stands for the contour, all rays making an angle of $\pi/4$ with the horizontal axis.
 $$ \nwarrow ~  \swarrow $$

\vspace{-.9cm}

$$ 0
$$

\vspace{-.9cm}

$$  \nearrow~ \searrow   $$
%
%
%
%
For $s=t$, the Pearcey kernel can also be written
\be
 K^{\PR}_{t,t}(x,y)
 =\frac{p(x)q''(y)-p'(x)q'(y)+p''(x)q(y)-tp(x)q(y)}{x-y}
 \label{Pkernel1}\ee
 with
 $$
 q(y):=\frac{i}{2\pi}\int_{-i\iy}^{i\iy}e^{-\frac{U^4}4+\frac{tU^2}2-Uy}dU
 ,
  \qquad
   p(x):=\frac{1}{2\pi i}\int_{X}e^{ \frac{V^4}4-\frac{tV^2}2+Vx}dV,
   $$
satisfying both, the differential equations (using integration by parts)
\be
 p'''(x)-tp'(x)+xp(x)=0  , \qquad q'''(y)-tq'(y)-yq(y)=0
 \label{ode},
 \ee
 and the heat equations
 \be
 \frac{\pl p}{\pl t}=-\frac{1}{2}p''(x),\quad \frac{\pl q}{\pl t}=\frac{1}{2}q''(y)
\label{heat}
,\ee
whereas $K^{\PR}_{t,t}$ satisfies the following equation
\be
\frac{\pl K^{\PR}_{t,t}}{\pl t}=\frac{1}{2}(-p'(x)q(y)+p(x)q'(y))
. \label{kernel-PDE}\ee
The latter follows from taking $\pl/\pl t$ of the kernel (\ref{Pkernel}), which has for effect to multiply the exponentials under the integral (\ref{Pkernel}) with $\frac 12(U^2-V^2)/(U-V)=\frac 12 (U+V)$.

 \vspace*{.9cm}

\begin{figure}
\includegraphics[width=100mm,height=120mm]{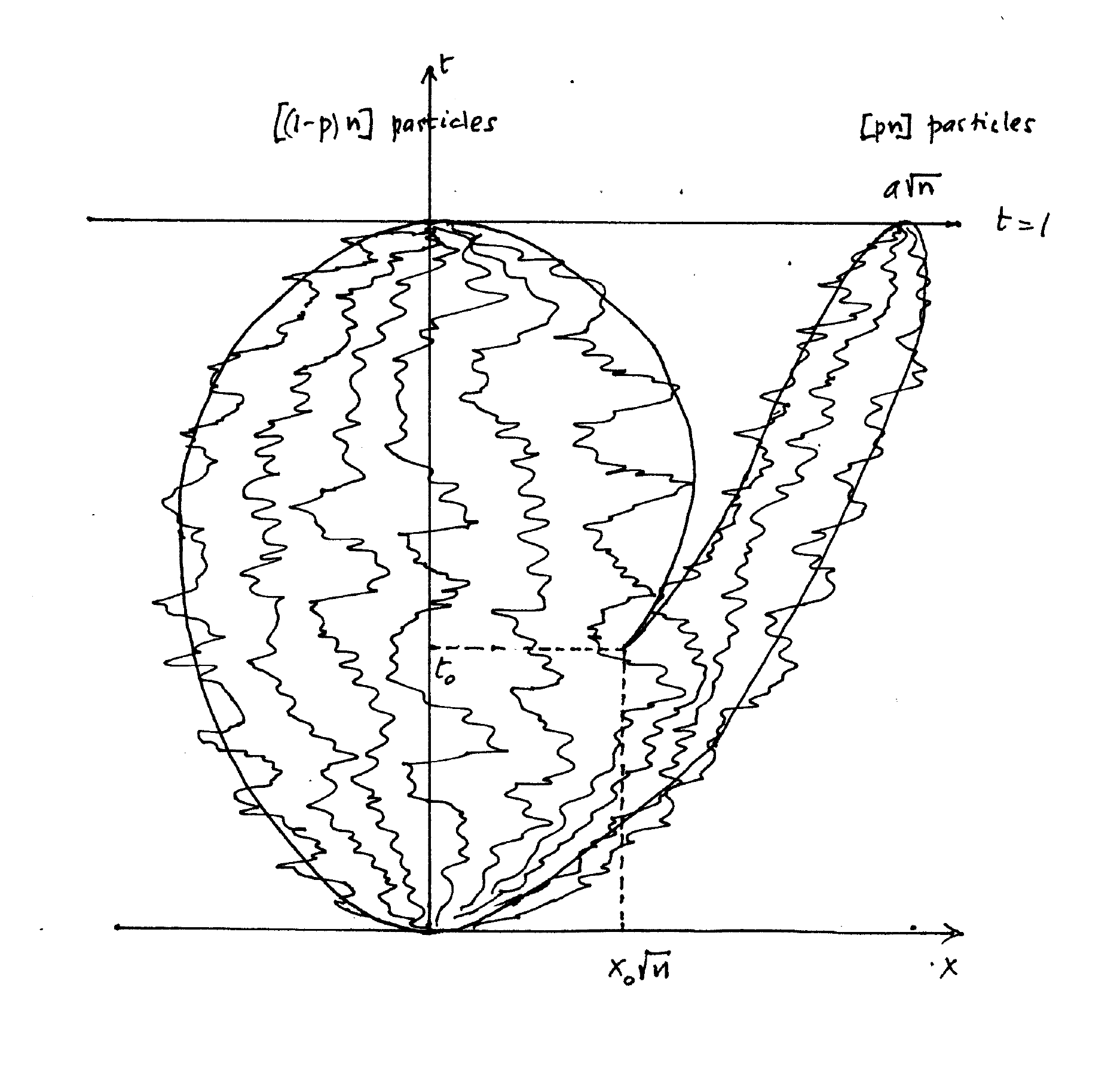}
\vspace*{-0.3cm}

\caption{The Pearcey process for $b=0$.}\label{fig1}\end{figure}

Consider $n$ non-intersecting Brownian motions, with
 $0<p<1$ and $b<a$:
%
%
\bean   
\BP^{(b,a)}_{n} \left(\!\!\!\!\begin{tabular}{c|c}
& all $x_j(0) =0$\\
  $\displaystyle \bigcap_{1\leq i\leq m}
  \left\{\mbox{all~}
 x_j \left(t_i\right)~ \mbox{for}~  1\leq j\leq n \right\} $
&
 $pn$ paths end up at $a$ at $t=1$\\
  &
$(1\!-\!p)n$ paths end up at $b$ at $t=1$
\end{tabular}\!\!\!\right)
\eean
When $n\rg \iy$, the mean density of Brownian particles
has its support on one interval for $t\sim  0$ and on
two intervals for $t\sim  1$, so that a bifurcation
appears for some intermediate time $t_0$, where one
interval splits into two intervals. At this point the
boundary of the support of the mean density has a cusp.
We show that near this cusp, the same Pearcey process
appears, independently of the values of $a$, $b$ and
$p$, showing ``{\em universality}" of the Pearcey process; see Figure \ref{fig1}. As it turns out, it
is convenient to introduce the parametrization
 \be
  p=\frac{1}{1+q^3} \mbox{~with~}  0< q< \iy  \mbox{~and let~}
  r:=\sqrt{q^2-q+1}.\label{notation}
   \ee

 \begin{theorem} \label{Theo 1}
 For $n\rg \iy$, the cloud of Brownian particles lie
within a region, having a cusp at location
$(x_0\sqrt{n},t_0)$, with
\be
  x_0 =\frac{(2a-b)q+(2b-a)}{q+1}t_0,~~~~
\frac{1}t_0= {1+2\left(\frac{r(a-b)}{q+1}\right)^2}.
 \label{x0,t0}\ee
  Moreover, the following probability tends to the
probability for the Pearcey process:
 \bea
\lefteqn{\lim_{n\rg\iy}\BP_n^{(b\sqrt{n},a\sqrt{n})}
 \left(\bigcap_{1\leq i\leq m}\left\{\mbox{all~}
x_j \left(t_0+ \left(\frac{c_0\mu}{n^{1/4}}\right)^2 {2\tau_i}
 \right)\in  {x_0}  n^{1/2}+c_0 A
\tau_i+\frac{c_0\mu}{n^{1/4}}{E^c}\right\} \right)
 }\no\\& &\hspace*{7cm}=\BP^{\cal P}\left(\bigcap_{1\leq i\leq m}\left\{{\cal P}(\tau_i ) \cap E =
\emptyset \right\}\right), \no\\\label{main}\eea
 using the following constants
 \bean
   \mu&=&\left(\frac{q^2-q+1
 }{q}\right)^{1/4}>0,~~~~~
 %
 %
 ,~~~~c_0:=\sqrt{\frac{ t_0(1-t_0)}{2}}=t_0\frac{r(a-b)}{q+1}>0.\eean
\be
 A = \frac{q^{1/2}(a-x_0)+q^{-1/2}(b-x_0)}{ (a-b)}
 \label{constants} \ee
\end{theorem}

In \cite{AvM-Pearcey}, Adler and van Moerbeke showed that the Pearcey transition probability (\ref{Pearcey tr-pr}) satisfies a non-linear PDE, expressible as a Wronskian of the expression (\ref{Pearcey PDE}) with some partial. This was obtained from taking a scaling limit, when $n\to \iy$, of the symmetric situation, i.e., where $b=-a$ and $p=1/2$. It came as a {\em surprise} to us that considering the asymmetric case leads to a different non-linear PDE, when $n\to \iy$, but nevertheless also expressible as a Wronskian of the same expression (\ref{Pearcey PDE}) with some other partial. A separate functional-theoretical argument then enables one to show that the expresssion (\ref{Pearcey PDE}) itself vanishes. This was one of the motivations for finding the exact scaling as presented in Theorem \ref{Theo 1}.

  To $E=\cup_{i=1}^r(y_{2i-1},y_{2i})\subset \BR$ one associates two operators, a divergence and an Euler operator
 $$\pl_{_{\!E}} =\sum_1^{2r}\frac{\pl}{\pl y_i},\qquad \vr_{_{\!E}} =\sum_1^{2r}y_i\frac{\pl}{\pl y_i}.$$




    \newpage

  \vspace*{0cm}

  \vspace*{-12cm}
  \bean
 \hspace*{7cm}
 &&\mbox{parametrization:}\\
 && q^3:=\frac{1-p}{p} ,~~ r:=\sqrt{q^2-q+1}
 .
   \eean

 \vspace*{1cm} \bean
  \hspace*{7cm} && ~~~~\mbox{cusp $x-x_0=2 \left(\frac{t-t_0}{3}\right)^{3/2}$ at}\\
  &&~x_0 
 =at_0\frac{2q-1}{q+1},~~
\frac{1}{t_0}= {1+2\left(\frac{ar}{q+1}\right)^2}.
 \eean

 \vspace*{5cm}





\begin{theorem}\label{Theo:2}

The log of the transition probability for the Pearcey process, which is non-stationary,
 \be
 {\mathbb Q}(t,E):= \log\BP^{\cal P}\left({\cal P}(t ) \cap E =
\emptyset 
 \right)
\label{Pearcey tr-pr}\ee
 satisfies the following 3rd order non-linear PDE in $t$ and the boundary points of $E$,
 \be
  \frac{\pl^3 {\mathbb Q}}{\pl t^3}
  +\frac{1}{8}\left(\vr_E-2t\frac{\pl}{\pl t}-2\right)\pl_{{_E}}^2 {\mathbb Q}
 -\frac{1}{2}\left\{\pl_{_E}^2 {\mathbb Q},\pl_{_E} \frac{\pl {\mathbb Q}}{\pl t}\right\}_{\pl_{_{\! E}}}
   =0
   , \label{Pearcey PDE}\ee
   with ``final condition", given by the Airy process\footnote{The Airy process is a stationary process, which describe the statistical fluctuations of the process about the curve appearing in Figure \ref{fig1}, away from the edge and properly rescaled. Its probability given at any time by the Tracy-Widom distribution ${\cal F} (x)$. The latter is given by the Fredholm determinant $\det(I-{\bf A})$ of the Airy kernel ${\bf A}$, restricted to the interval under consideration.} which is a stationary process, (by moving far out along the cusp
    $x=2 \left(\frac{t}{3}\right)^{3/2}$) 
%
   %
%
   $$\lim_{t\rg \iy} \BP^{\cal P}\left(\frac{{\cal P}(t)-2 \left(\frac{t}{3}\right)^{3/2}}{(3t)^{1/6}}\cap (-E)=\emptyset
   \right)=\det(I-{\bf A})_{(-E)}
   $$

 \end{theorem}

\remark It is interesting to compare the Pearcey PDE with the Airy process PDE; namely for semi-infinite
  intervals $E_1$ and $E_2$, the 3rd order non-linear PDE for the Airy
  joint probability
$$
 {\mathbb Q}^{\cal A}(t;x,y):=
  \log \BP^{\cal A}\left({\cal A}(t_1) \leq\frac{y+x}{2}, ~{\cal A}(t_2) \leq
  \frac{y-x}{2}\right),\mbox{   for  } t=t_2-t_1,
 $$
%
 reads
  \be
  2t   \frac{\pl^3 {\mathbb Q}^{\cal A}}{\pl t \pl x \pl y}
 =
\left(t^2 \frac{\pl}{\pl x}- x  \frac{\pl}{\pl y}\right)
 \left(\frac{\pl^2 {\mathbb Q}^{\cal A}}{\pl x^2}-
  \frac{\pl^2 {\mathbb Q}^{\cal A}}{\pl y^2}\right)
+8 \left\{  \frac{\pl^2 {\mathbb Q}^{\cal A}}{\pl x \pl y} ,
 \frac{\pl^2 {\mathbb Q}^{\cal A}}{\pl y^2} \right\}_y,
\label{Airy PDE1}\ee
with "final condition":
$$\lim_{t_2-t_1\to \iy}
 \BP^{\cal A}\left({\cal A}(t_1) \leq u_1, ~{\cal A}(t_2) \leq
 u_2\right)={\cal F}(u_1)  {\cal F}(u_2)
 . $$






  
  In the last section (section \ref{section 7}), we develop -in a formal way- the central role played by the spectral curve (or Pastur equation \cite{Pastur}) in the steepest
descent analysis used to prove the universal behavior of the kernel as $N \to \infty$ for the different problems of
non-intersecting Brownian motions. The spectral curve is precisely the function which appears in the steepest descent analysis. The spectral curve associated
to the problem provides the universal limiting kernel obtained after a proper rescaling of the variable
around a singularity of the problem.

\section{Non-intersecting Brownian motions on $\BR$, forced to several points} \label{section1}
In the expression below,
 ${\cal H}_{n}(E)$ is the set of all Hermitian
matrices with all eigenvalues in $E$. Note that in
general one has the following, using the Karlin-McGregor formula\footnote{$\Dt_n(x_1,\ldots,x_n)$ is the Vandermonde determinant.} (see \cite{Karlin,Brezin2,Brezin3,Brezin4,Brezin5,TW-Pearcey,BleKui2}):
\bea  \lefteqn{
\BP^{(b_1,\ldots,b_p)}_{n} \left(\begin{tabular}{c|c}
& all $x_j(0) =0$\\
all $x_j(t)\in E$ for $1\leq j\leq n$ &
 $n_1$ paths end up at $b_1$ at $t=1$\\
  & \vdots\\
  &
$n_p$ paths end up at $b_p$ at $t=1$
\end{tabular}\right)}
\no\\
 &=&\!\!\!\!\!\!\!\!\! \lim_{\tiny\begin{array}{c}
  \mbox{\tiny   all $~\gamma_i\rightarrow 0$}
  \\
 \mbox{\tiny $\delta_1, \ldots , \delta_{n_1}\rightarrow b_1$}
  \\\vdots\\
  \mbox{\tiny $\delta_{n_1+\ldots+n_{p-1}+1}, \ldots ,
  \delta_{n}
  \rightarrow b_p$}
  \end{array}}\!\!\!\!\!\!\!\!\!   {\displaystyle\int_{E^n}
  \frac{\prod_1^n dx_i}{Z_{n}}}
  %
  \det\left(p(t;\gamma_i,x_j)\right)_{1\leq i,j\leq
n}\det(p(1\!-\!t;x_{i'},\delta_{j'}))_{1\leq i',j'\leq
n}
  ,  \no\\
  %
%
 &=&\frac{1}{Z_n}\left.\int_{\tilde E^n}
 \Delta_n(x_1,...,x_n)
%
 \prod^p_{\ell=1} \Delta_{n_{\ell}}
(x^{(\ell)})\prod_{j=1}^{n_{\ell} }
e^{-\frac{1}{2}x_j^{(\ell)^2}+\tilde
b_{\ell}x_j^{(\ell)}
 }dx_j^{(\ell)} \right|
 _{\begin{array}{c}\tilde
 E=E\sqrt{\frac{2}{ t(1-t) }}
    \\
    \tilde b_\ell=\sqrt{\frac{2t}{1-t}}b_\ell
    \end{array}
}
 \no\\ \no\\&=&
 \frac{1}{Z_n}\int_{\HR_n\left(E\sqrt{\frac{2}{t(1-t)}}\right)}
 dM e^{-\frac{1}{2}\Tr (M^2-2A_t M)}dM
 ~~~ \no\\\no \\\no \\
 &=&
 \frac{1}{Z_n}\det\left(
       \begin{array}{l}
          \left({\displaystyle\int_{\tilde E}}~x^{i+j}
          e^{-\frac{x^2}{2}+
    \tilde b_{1}x}dx\right)_{0\leq i \leq n_1-1,
          ~0\leq j\leq n-1 }\\
          \hspace{1cm} \vdots\\
          \left({\displaystyle\int_{\tilde E}}~x^{i+j}e^{-\frac{x^2}{2}+
   \tilde b_{p}x}dx\right)_{0\leq i \leq n_p-1,
          ~0\leq j\leq n-1 }\\
       \end{array}
     \right) \no\\ \no\\
&=& ~~~ { {\det(I-H_n^{(p)})}_{ E^c },~~~(\mbox{\bf
Fredholm determinant})}
  \label{Brownian transition}\eea
%
where $H_n^{(p)}(x,y)$ is the kernel  (setting
$t_k=t_\ell=t $)
 \bea
\lefteqn{H_n^{(p)}(x,y)dy }\label{BaikKernel}
   \no\\
&=&
 -\frac{dy}{2 \pi^2 \sqrt{(1-t_k)(1-t_\ell)}}
  \int_{\cal C} dV
\int_{L+i\BR} dU ~\frac{e^{ -\frac{ t_kV^2}{1-t_k}  +
\frac{2xV}{1-t_k}}}{e^{ - \frac{t_\ell U^2}{1-t_\ell} +
\frac{2yU}{1-t_\ell} }}
 \prod_{r=1}^p\bigg(\frac{U-b_r}{V-b_r}\bigg)^{n_r}
\frac{1}{U-V}
\no\\ \no\\
 &&-\left\{ \begin{array}{l}
            0 ~~~~~~~~~~~~~~~~~~~~~~~~~~~~~~~~~~~~
            \mbox{for}~~t_k\geq
            t_\ell\\  \\  
             \frac{1}{\sqrt{\pi (t_\ell-t_k)}}
            e^{-\frac{(x-y)^2}{t_\ell-t_k}}
            e^{\frac{x^2}{1-t_k}-\frac{y^2}{1-t_\ell}},~~~~~
            \mbox{for}~~t_k<t_\ell
            \end{array}\right..  \label{BM kernel}\eea
 where $X$ is a contour consisting of the two incoming rays from $\pm \iy e^{i\pi/4}$ to $0$ and the two outgoing rays from $0$ to $\pm \iy e^{-i\pi/4}$, provided no $b_r=0$.
  In the expression above, $A_t$ is the diagonal matrix \be
 A_t:= \left(\begin{array}{cccccccccccc}
   \tilde b_1\\
   &\ddots& & & &{\bf O}\\
   & &\tilde b_1\\
   & & &\tilde b_2\\
   &{\bf O}& & &\ddots\\
   & & & & &\tilde b_2&                 \\
   & & & & &          &\ddots&          \\
   & & & & &          &       &\tilde b_{p}&\\
   & & & & &          &       &            &\ddots&\\
   & & & & &          &       &            &      &\tilde b_p\\
\end{array}\right)\begin{array}{l}
\updownarrow  n_1\\
\\
\\
\updownarrow n_2
 \\
 \vdots
\\ \\
\\
\updownarrow n_p
\end{array}
  \label{diagmatrix}~  \mbox{with}~~~
   \tilde b_i=b_i \sqrt{\frac{2t }{1-t}},\ee


The main expression appearing in (\ref{Brownian transition}) contains the matrix integral
  \bea
  \BP_n(E;b_1,\ldots,b_p)&=&\frac{1}{Z_n} \int_{ E^n}
   \Delta_n(x_1,...,x_n)
%
 \prod^p_{\ell=1} \Delta_{n_{\ell}}
(x^{(\ell)})\prod_{j=1}^{n_{\ell} }
e^{-\frac{1}{2}x_j^{(\ell)^2}+ b_{\ell}x_j^{(\ell)}
 }dx_j^{(\ell)}
%
 \no\\ \no\\&=&
 \frac{1}{Z_n}\int_{\HR_n\left(E\right)}
 dM e^{-\frac{1}{2}\Tr (M^2-2A M)}
 , \label{integral}\eea
 %
  which is now viewed as a function of the boundary points of $E$
 and the target points $b_i$, for which one assumes a linear
 dependence $${\sum_1^{p}}c_ib_i=0\mbox{~~with~~}{\sum_1^{p}}c_i=1.$$
Introduce the following operators:
 \bea
  \pl_{_{\!E}} &:=&
\left\{\begin{array}{l}\mbox{sum of partials
in the}\\
\mbox{boundary points of}~ E\end{array}\right\}
 \no\\
 \varepsilon &:=&\left\{\begin{array}{l}\mbox{Euler operator
in the}\\ \mbox{boundary points of}~ E
\end{array}\right\}-\sum_1^{p-1} b_i\frac{\pl}{\pl
b_i} \no\\
\pl^{(\ell)}_{b}&:=& c_\ell\sum_1^{p-1}\frac{\pl}{\pl
b_i} -\frac{\pl}{\pl b_\ell}(1-\delta_{\ell p})
 ,~~~~\mbox{one checks}~\sum_{\ell=1}^p\pl^{(\ell)}_{b}=0  
\label{del} \eea 

{\bf Proposition:} \cite{AvM-Mops} The expression $\log \BP_n$ satisfies a
non-linear PDE in the boundary points of the interval
$E$ and in the target points $b_i$, given by the ({\em near-Wronskian})
determinant of a $(p+1)\times (p+1)$ matrix  {
\be \det \left(
  \begin{array}{lllllllllll}
  F_1&F_2&F_3&\ldots&F_p&0 \\
  F_1' &F_2'&F_3'& \ldots&F_p'&G_1 \\
    F_1'' &F_2''&F_3''& \ldots&F_p''&G_2 \\
    \vdots&\vdots&\vdots&&\vdots&\vdots
    \\
    F_1^{(p )} &F_2^{(p )}&F_3^{(p )}& \ldots&F_p^{(p )}&G_p \\
    %
    %
    %
\end{array}
\right)=0 ,~~~~~{}':=\pl_{_{\!E}}  \label{7.4.11}
 \ee}
 where the $F_\ell$ and $G_\ell$ are given by\footnote{ with $C_\ell=-2n_\ell\left((1-c_\ell)b_\ell+\sum_{j\neq \ell}\frac{n_j}{b_\ell-b_j}\right).$}  ($G_0=0$)
%
%
\bea F_{\ell} &:=&
 \left(  \pl^{(\ell)}_{b} +c_\ell\pl_{_{\!E}}\right)\pl_{_{\!E}} \log
\BP_n+n_{\ell}
\no\\
G_{\ell+1}&:=&\pl_{_{\!E}} G_\ell
   {~~+\sum_{j=1}^p
  {
(\pl_{_{\!E}}^\ell F_{j})} \left(
  \pl_{_{\!E}}\frac{H_{ j }^{(1)}}{F_{j}}-
\pl^{(\ell)}_{b}\frac{H_{ j }^{(2)}}{F_{j}}
%
\right)},
\no\\
 H_{ \ell }^{(1)} &:=&
 \left(-c_\ell\pl_{_{\!E}}
 \varepsilon+((\varepsilon-1)c_\ell+2) (\pl^{(\ell)}_{b}
 +c_\ell\pl_{_{\!E}} )\right)\log \BP_n
  +C_{ \ell }
 \no\\
H_{ \ell }^{(2)} &:=& \left(1-\varepsilon
 +2b_\ell \pl_{_{\!E}}
\right)\left(\pl^{(\ell)}_{b} +c_\ell\pl_{_{\!E}}\right)
\log\BP_n .
 \label{F,G,H} \eea

  \noindent{\em Example:} Setting ${H}_{\ell}:=
 \left\{H_{ \ell }^{(1)},F_{\ell}\right\}_{\pl_{_{\!E}}}-
\left\{H_{ \ell }^{(2)}, F_{\ell}\right\}_{
\pl^{(b)}_{\ell}}$, one has:
\begin{eqnarray}
 p=1
&\Longrightarrow& \det\left(\begin{array}{cc}
   F_1 & 0\\
   F_1'& \frac{H_1}{F_1}\end{array}\right)
=H_1=0
\label{eqts}\\
p=2 &\Longrightarrow& F_1F_2\det\left(\begin{array}{ccc}
   F_1 & F_2&0\\
   F_1'&F'_2& \frac{H_1}{F_1}+\frac{H_2}{F_2}\\
  F_1''&F''_2& \frac{H'_1}{F_1}+\frac{H'_2}{F_2}
  \end{array}\right)=
 (H_1F_2+H_2F_1)\{F_1,F_2\}'
 \no\\
 &&~~~~~~~~~~~~~~~~~~~~~~~~~~~~~~~~~~~~~~~~~
 -
  (H'_1F_2+H'_2F_1)\{F_1,F_2\}=0\no
\end{eqnarray}
%


\section{Proof of Theorem \ref{Theo 1}} \label{section2}

\underline{Proof}: It is easily computed by the
Pastur-Marcenko method \cite{Pastur}, which states that,
given a diagonal matrix $A=(a_1,\ldots,a_n)$ and its
 spectral function $d\sigma(\lb):=\frac{1}{n}\sum_i \dt(\lb-a_i)$,
 the random Hermitian ensemble with probability defined by
\be
 \frac{1}{Z_n}\int_{\HR_n(E )}dM e^{-\frac{n}{2v^2}\Tr
(M-A)^2}
 ,\label{Pastur}\ee
 has, in the limit when $n\rg \iy$, a spectral density $d\nu(\lb)$, whose Stieltjes
 transform
  $$
   f(z)=\int_{-\iy}^{\iy} \frac{d\nu(\lb)}{\lb -z},~~\Im m
   ~z\neq 0,
    $$
     satisfies the integral equation
      \be
      f(z)=\int_{-\iy}^{\iy}
      \frac{d\sigma (\lb)}{\lb-z-v^2 f(z)}
      ,\label{IntEqt}\ee
    and, in view of the spectral function $d\sigma(\lb):=\frac{1}{n}\sum_i \dt(\lb-a_i)$, this becomes,
    \be
    f(z)=\sum_{i=1}^n \frac{1/n}{a_i-z-v^2 f(z)}
   .\label{AlgEqt}\ee
    Consider the situation (\ref{Brownian transition}), where $n_i$ particles are forced to $b_i$. Then, upon setting the variance $v^2=1$ and defining $g(z):=f(z)+z$, the equation (\ref{AlgEqt}) reads\footnote{See Section 7, formula (\ref{fraction numbers}), for further comments on the ``fraction numbers" $\vr_i$.}
    \be
    g-z+\sum_{i=1}^p \frac{\vr_i}{g-\tilde b_i }=0,~~\mbox{with}~\vr_i=\frac{n_i}n.
    \label{AlgEqt1}\ee
    The density of the equilibrium
    distribution is then given by
     $$
      \frac{d\nu(z)}{dz}=\frac{1 }{\pi} |\Im m~ f(z)|
       = \frac{1 }{\pi} |\Im m~ g(z)| ~, \mbox{ ~~  with  } z\in \BR
     .$$
     %
     %
   {\em Remark}: It is precisely this Pastur-Marchenko equation (\ref{AlgEqt1}), which will appear in the argument of steepest descent for the corresponding kernel in section 7. Namely, equation (\ref{AlgEqt1}) is the derivative (\ref{S'}) of the S-function, appearing in (\ref{S}).

     We will now specialize to two target points $p=2$, where $n_1=pn$ and $n_2=(1-p)n$.
 From (\ref{Brownian transition}), it follows that for non-intersecting Brownian motions forced to $b\sqrt{n}<a\sqrt{n}$, one has
\be
 \BP_n^{b\sqrt{n},a\sqrt{n}} \left( \mbox{all $x_j(t)\in \sqrt{n}E$} \right)
 =\frac{1}{Z_n} \int_{{\cal
H}_{n}\left(E\sqrt{\frac{2}{t(1-t)}}\right)}
   e^{-\frac{n}{2}\Tr( M-A_t)^2 }dM,
\label{Brownian transition1}\ee
with
\be
A_t
 := \left(\begin{array}{cccccc}
   a \sqrt{\frac{2t }{1-t}}=\a\\
   &\ddots& & & &{\bf O}\\
   & &a \sqrt{\frac{2t }{1-t}}=\a\\
   & & &b \sqrt{\frac{2t }{1-t}}=\b\\
   &{\bf O}& & &\ddots\\
   & & & & &b \sqrt{\frac{2t }{1-t}}=\b
\end{array}\right)\begin{array}{l}
\updownarrow pn\\
\\
\\
\\
\updownarrow (1-p)n
\end{array},
  \label{0.1}\ee
 For $A_t$ as above, the integral equation (\ref{IntEqt})
     becomes
  the algebraic equation 
    $$g-z+{\frac {1-p}{g-\b}}+{\frac
{p}{g-\a}}=0
 $$
%
which, upon clearing, leads to a cubic equation for
$g$,
 \be
 G(g):={g}^{3}- \left(  z+\a+\b \right)
 {g}^{2}+ \left( z(\a+\b
 )+\a\b +1 \right) g-(\a \b z+ (1-p)\a
  +p\b)=0,
\label{G-equation}\ee
with roots given by
 $
 g= \tilde q+\sqrt[3]{\tilde r+\sqrt{\Dt_1}}+\sqrt[3]{\tilde r-\sqrt{\Dt_1}}
 $,
 with a quartic discriminant in $z$,
$$\Dt_1(z)=(\a-\b)^2\prod_0^3(z-z_i)=0;$$
$\tilde q,\tilde r$ are polynomials of $z,\a, \b$. Thus one finds the
following
      $$
         \frac{d\nu(z)}{dz}=\frac{1}{\pi}\left| \Im m~ g(z) \right|
        = \left\{\begin{array}{l}
        {\displaystyle \frac{1}{\pi}|\Im m~ g(z)|\mbox{   for $z$ such that }
        \Dt_1(z) <0
         }\\ \\
         0 \mbox{   for $z$ such that }
        \Dt_1(z) \geq 0
         \end{array}
         \right.
         $$
Therefore the support of the equilibrium measure will be
given by
 \bean
 &&\mbox{either two intervals}~[z_2,z_0]\cup [z_1,z_3]\\
 &&\mbox{or two intervals touching}~[z_2,z_0]\cup [z_0,z_3]\\
 &&\mbox{or one interval}~[z_2,z_3]
 \eean
 for the real roots of
$\Dt_1(z)=(\a-\b)^2\prod_0^3(z-z_i)=0$. Thus depending
on the values of the parameters $\a,~\b$ and $p$, there
will be four real roots or two real roots, with a
critical situation where two of the four real ones
collide, say $z_1=z_0$. The latter occurs exactly when
the discriminant $\Dt_2$ (with regard to $z$) of $\Dt_1$
vanishes, namely when
$$\Dt_2(\a,\b,p)=4 p(1-p)\rho
 \Bigr((\rho-1)^3-27p(1-p)\rho\Bigl)^3\Bigr|_{\rho=(\a-\b)^2}
 =0.$$
 This polynomial has one positive root
 (the others being imaginary), and one checks that,
 taking into account $\a>\b$, and defining $p$  and $r$ as in (\ref{notation}),
  \be
   \a-\b=\sqrt{\rho}=
    \Bigl(3 {p}^{1/3} \left( 1\!-\!p \right) ^{2/3}+3
   {p}^{2/3}({1\!-\!p})^{1/3}+1\Bigr)^{1/2}= \frac{q\!+\!1}{r}
  >0 
   \label{bifurc condition}\ee 
  For this precise value of the parameter $\a-\b$, two of the four
  boundary points of the support coincide, namely the
  roots $z_0$ and $z_1$ of $\Dt_1(z)$ coincide:
  \be
z_0=z_1=\b+ \frac {2q-1}{r}
 =\a+ \frac {q-2}{r}
.\label{double root}\ee
This double root $z_0$ is found by stating that, under the
condition $\a -\b=(q+1)/r$, the
polynomials $\Dt_1$ and $\Dt_1'$ have a common root; in
other terms for appropriate choices of
$c_1,c_2,c_3$, some linear combination of the polynomials $\Dt_1$ and $\Dt'_1$ becomes a linear polynomial
 $
  (4z+c_1)\Dt_1(z)-(z^2+c_2z+c_3)\Dt'_1(z)=
    c_4z+z_5
    $
    for some $c_4$ and $c_5$. Since $z_0$ is a root of the left hand side, it also must be a root of the right hand side. This is to say $z_0=-c_5/c_4$, yielding the expression (\ref{double root}).

The critical time $t_0$ is then obtained from setting $t=t_0$ in (\ref{0.1}),
\be
 \a= a
\sqrt{\frac{2t_0 }{1-t_0}} \mbox{~~~and~~~}\b=b \sqrt{\frac{2t_0
}{1-t_0}}
 \label{alpha and beta}\ee
  from which one computes $t_0$, by taking the difference and by using (\ref{bifurc condition}),
$$t_0=\frac{(q+1)^2}{(q+1)^2+2(a-b)^2(q^2-q+1)},$$
and from which one further computes 
\be
 c_0:=\sqrt{\frac { t_0(1-t_0)}{2}}=t_0\frac{r(a-b)}{q+1} ~~\mbox{and}~~\frac{t_0}{c_0}=\sqrt{\frac{2t_0}{1-t_0}}
   =\frac{q+1}{(a-b)r}, \label{c0}\ee
confirming the expression for $c_0$ in (\ref{constants}).
Then from (\ref{alpha and beta}) and (\ref{c0}) one deduces
\be
 \a=a \frac{t_0}{c_0}= \frac{a(q+1)}{(a-b)r}~~~~~~\b=b \frac{t_0}{c_0}= \frac{b(q+1)}{(a-b)r}
 \label {alpha-beta}  .\ee
Defining $x_0:=z_0c_0$, one computes from (\ref{double root}),
 \be
  z_0=\frac{x_0}{c_0}=\frac{(2a-b)q+(2b-a)}{(a-b)r
  }
 .\label{z0}\ee
 Next, one computes the double root of the $G$-equation (\ref{G-equation}) for the value $z=z_0$. Indeed, using (\ref{double root}) and (\ref{alpha-beta}), one checks that for some root $g=g_0$,
  $$
 G= \frac{\pl G}{\pl g}=\frac{\pl^2G}{\pl g^2}  = 0,
 $$
 where
  \be  g_0:=\frac{1}{3}(z_0+\alpha+\beta)=
   \frac{1}{3}(2a\frac{t_0}{c_0}+b\frac{t_0}{c_0}+\frac{q-2}{r})=\frac{aq+b}{(a-b)r}
 ,\label{g-root}\ee
 and from (\ref{alpha-beta}) and $b<a$ that
  \be \beta<g_0<\a \label{root2}.   \ee
The point $z_0$ in (\ref{z0}) refers to the matrix integral variables on the right hand side of (\ref{Brownian transition1}), which can then be transformed
into the Brownian motion variables, according to $ E=\tilde E\sqrt{\frac{ t(1-t) }{2}}$ at $t=t_0$, which gives the transformation from the matrix integral variables to the Brownian motion variables. 
      Hence the critical point in the Brownian picture takes place at (using (\ref{c0}) and (\ref{z0}))
   $$  (x_0\sqrt{n},t_0)=(z_0c_0\sqrt{n},t_0)=
    \left(\frac{(2a-b)q+(2b-a)}{q+1}t_0\sqrt{n},~ t_0\right)
   .$$


To prove the main statement (\ref{main}) in Theorem \ref{Theo 1}, start with the
kernel (\ref{BM kernel}) of the non-intersecting
Brownian motion:
 \bea
\lefteqn{H_n(x,y;t_k,t_\ell)dy}\no\\
 &=&-\frac{dy}{2\pi^2\sqrt{(1-t_k)(1-t_\ell)}}
 \int_{\CR}dV\int_{L+i\BR}dU~\frac{e^{-\frac{tV^2}{1-t_k}+
\frac{2xV}{1-t_k}}}{e^{
-\frac{tU^2}{1-t_\ell}+\frac{2yU}{1-t_\ell}} } \no \\
& &\qquad\qquad\qquad\qquad\qquad\times\left(\frac{U-b}{V-b}\right)^{n_2}
\left(\frac{U-a}{V-a}\right)^{n_1}\frac{1}{U-V} . \label{BM-kernel}\eea
One first needs to prove that for some $\varphi_n(\lb,\tau)$:
  \bean\lefteqn{\lim_{n\rg \iy}
 \varphi_n(x,t_k)H_n(x,y;t_k,t_\ell)\varphi_n(y,t_\ell)^{-1}\Bigr|_{\mbox{\tiny rescaling}}
 }\\
 &=&K^{\PR}(\xi,\eta; \tau_k,\tau_\ell)\\
 & :=&
 - \frac{1 }{4\pi^2 }\int_X d\om_v
\int^{i\iy}_{-i\iy}d\om_u ~~ 
\frac{e^{-\frac{\om_u^4}{4}+\frac{\tau_\ell\om_u^2}{2}-\om_u \eta}}
 {e^{-\frac{\om_v^4}{4}+\frac{\tau_k\om_v^2}{2}-\om_v\xi}}
~\frac{1}{\om_u-\om_v}
\eean
with the rescaling $\left\{ x\atop y\right\}$
\beq\label{rescaling1}
  t_i= t_0+ (c_0\mu)^2\frac{2\tau_i}
{n^{1/2}}
,~~
\left\{ x\atop y\right\}= c_0\left( {z_0}  n^{1/2}+ A
\left\{ \tau_k\atop \tau_\ell\right\}+\mu\frac{\left\{ \xi\atop \eta\right\}}{n^{1/4}}\right),~~
\eeq
 with constants $A, ~\mu,$ given by
    (\ref{constants}).

    Consider the change of variables
$
  U:= \frac{c_0u\sqrt{n}}{t_0}.
 $
The form of the Brownian motion kernel (\ref{BM kernel}) suggests, by putting the two $U$-factors of the integrand in the exponential, the function $F(u)$, which one observes, at leading order, is closely related to the function $G(u)$, defined in (\ref{G-equation}); namely\footnote{The function $F(u)$ should contain the term $\log \frac{\sqrt{n}c_0}{t_0}$; however as the same rescaling is made in the $v$-variables, this same term will appear with a different sign and therefore they will cancel. Consequently, this term will be omitted.}
 \bea
F(u) &:=&\frac{u^2}{2}-uz
+p\log(u-\a)
+(1\!-\!p)\log(u- \beta)
%
\Bigr|_{\a=\frac{at_0}{c_0},\quad
                            \b=\frac{bt_0}{c_0},\quad
                             z=\frac{x_0}{c_0}},\label{F-function}\eea
  with
\bean
  F'(u)=\left.\frac{G(u)}{(u-\a)(u-\b)}
 \right|_{
                           \a=\frac{at_0}{c_0},\quad
                            \b=\frac{bt_0}{c_0},\quad
                             z=\frac{x_0}{c_0}
                              }
. \eean
 Remember from (\ref{g-root}) that
 \be u_0:= g_0
 = \frac{aq+b}{(a-b)r}\label{u0}\ee
is a root of $G(u)=0$ and two of its derivatives,
 $$ G(u_0)=G'(u_0)=G''(u_0)=0\mbox{  and  }
 G'''(u_0)=6$$
and then one computes, since\footnote{upon using (\ref{double root}), $\a -\b=(q+1)/r$ and the root $u_0=\frac{1}{3}(z_0+\a +\b)$ as in (\ref{g-root}).} $(u_0-\a)(u_0-\b)=-q(q^2-q+1)^{-1}=-\mu^{-4}$,
\beq \label{derivFnul}
  F'(u_0)=F''(u_0)=F'''(u_0)=0  \mbox{  and  }
  \frac{1}{4!}F^{(iv)}(u_0)=   -\frac{(q^2-q+1)}{4q}
\eeq
   and so
 \bean
F(u)
&=&F(u_0)-\frac{\mu^4}{4}\left( {u-u_0} \right)^4+{\bf
O}(u-u_0)^5.
\eean
In the calculation below, the first equality $\stackrel{*}{=}$ is obtained by
doing all the substitions below, except for the last one
$u\mapsto \om_u$, whereas the second equality
$\stackrel{**}{=}$ is obtained by the substitution to
the new integration variable $u\mapsto \om_u$; the
expression after $\stackrel{**}{=}$ contains a term
$n^{1/4} \om_u$, which contains the new integration
variable and which blows up as $n^{1/4}$ . Hence this
coefficient must be put $=0$, yielding the value of $A$
as in (\ref{constants}). In the next equality one uses this
value, thus yielding in the end,
     \bean
\lefteqn{\left.\begin{array}{l}\displaystyle\frac{tU^2}{1\!-\!t}-\frac{2Uy}{1\!-\!t}
+n_2\log(U\!-\!b)+n_1\log(U\!-\!a)\\ \\ \\ \end{array}
\right|_{\begin{array}{l}
     n_1=np,~~n_2=(1-p)n\\
     U= \frac{c_0u\sqrt{n}}{t_0}\\
     t=t_0+(c_0\mu)^2\frac{2\tau_\ell}{n^{1/2}}\\
     a\mapsto a\sqrt{n},~~
     b\mapsto b\sqrt{n}\\
     y= c_0\left( {z_0}  n^{1/2}\!+\! A
\tau_\ell\!+\!\mu\frac{\eta}{n^{1/4}}\right)
     \\
     u=u_0+\frac{ \om_u}{\mu n^{1/4}}\end{array}}
     } \\
%
 \\
 &\stackrel{*}{=}&nF(u)+{n}^{1/2}~ {\tau_\ell u\mu^2}
 \left(\frac{u}{2}-\frac{t_0}{c_0}x_0-\frac{A}{\mu^2}\right)
   -n^{1/4}{u\eta   \mu}
   \\
   &&\hspace*{3cm}+
  {t_0 u\mu^4} \tau_\ell^2\left(\frac{u}{2}-\frac{t_0}{c_0}x_0-\frac{A}{\mu^2}\right)
 +O(n^{-1/4}) \\
 &\stackrel{**}{=}&
 nF(u_0)+{n}^{1/2}~ {\tau_\ell u_0\mu^2}
  \left(\frac{u_0}{2}-\frac{t_0}{c_0}x_0-\frac{A}{\mu^2}\right)\\
&& \hspace*{3cm}+
  n^{1/4}
 \bigl( {\tau_\ell}  \om_u( u_0-
\frac{t_0}{c_0}x_0-\frac{A}{\mu^2})-\eta u_0\bigr)\mu+O(1)
\\
 &=&
 nF(u_0)+{n}^{1/2}   \frac{\tau_\ell
u_0^2\mu^2}{2}-
  n^{1/4}
 { u_0\mu \eta}
-
{{ \left( \frac{\om_u^4}{4}- \frac{\tau_\ell\om_u^2 }{2 } +
\eta \om_u \right)} }-
\frac{t_0u_0^2\mu^4}{2}\tau_\ell^2
\\
&& +O(n^{-1/4}),\eean
%
  using the value (\ref{constants}) of
  $$
 A=     \sqrt{q}
 +\frac{t_0}{c_0}\mu^{2}\left( {b} -  x_0\right)
 = \frac{q^{1/2}(a-x_0)+q^{-1/2}(b-x_0)}{ (a-b)}
.$$ Similarly,
\bean
{\left.\begin{array}{l}\displaystyle\frac{tV^2}{1\!-\!t}
-\frac{2Vx}{1\!-\!t} +n_2\log(V\!-\!b)+n_1\log(V\!-\!a)\\ \\ \\
\end{array} \right|_{\begin{array}{l}
     n_1=np,~~n_2=(1-p)n\\
     V= \frac{c_0v\sqrt{n}}{t_0} \\
     t=t_0+(c_0\mu)^2\frac{2\tau_k}{n^{1/2}}\\
     a\mapsto a\sqrt{n},~~
     b\mapsto b\sqrt{n}\\
     x= c_0\left( {z_0}  n^{1/2}\!+\! A
\tau_k+\mu\frac{\xi}{n^{1/4}}\right)
     \\
     v=v_0+\frac{ \om_v}{\mu n^{1/4}}\end{array}}
     }\eean\bean
 &=&
 nF(u_0)+{n}^{1/2}   \frac{\tau_k
u_0^2\mu^2}{2}-
  n^{1/4}
 { u_0\mu \xi}
-
 { \left(\frac{\om_v^4}{4}- \frac{\tau_k\om_v^2~}{2~}+
\xi \om_v \right)}  -
\frac{t_0u_0^2\mu^4}{2}\tau_k^2
 \\
 &&+O(n^{-1/4})\eean
Moreover, using $\frac{2c_0^2 }{ (t_0-1)t_0}=-1$, together with the rescalings above, one finds
$$
 \frac{{ dy}}{2\pi^2\sqrt{(1-t_k)(1-t_\ell)}}\frac{dUdV}{U-V}
  =-\frac{ {  d\eta}}{4\pi^2}
  \frac{d\om_ud\om_v}{ (\om_u-\om_v)}+O(n^{-1/2})
  $$
Summarizing, one obtains:
  \bean
 \lefteqn{ \hspace*{-.2cm}\frac{-{dy}}{2\pi^2 \sqrt{(1-t_k)(1-t_\ell)}
 }
 \int_{\cal C}dV\int_{L+i\BR}dU~
 \frac{e^{-\frac{t_kV^2}{1-t_k}+
\frac{2xV}{1-t_k}}}
{e^{-
\frac{t_\ell U^2}{1-t_\ell}+\frac{2yU}{1-t_\ell}}
}}\\
 && \qquad\qquad\qquad\qquad\qquad\qquad\qquad
 \times \left(\frac{U-b}{V-b}\right)^{n_2}
\left(\frac{U-a}{V-a}\right)^{n_1}\frac{1}{U-V} \\
&=&
   e^{
 { u_0\mu( \xi-\eta)} n^{1/4}}
     %
     e^{\frac{1}{2}\sqrt{n}u_0^2\mu^2(\tau_k -\tau_\ell) }
     %
     e^{\frac{1}{2}t_0u_0^2\mu^4(\tau_k^2-\tau_\ell^2)}
     %
  \\
  &&\frac{{d\eta}}{4\pi^2}
\int_X d\om_v \int_{-i\iy}^{i\iy} d\om_u
  e^{ \frac{\om_v^4}{4}- \frac{\tau_k\om_v^2~}{2~}+
\xi \om_v -\frac{\om_u^4}{4}+ \frac{\tau_\ell\om_u^2~}{2~}-
\eta \om_u } \frac{1}{\om_u-\om_v}+O(n^{-1/4})
 ,\eean
  where we deformed the $u,v$ contours, by translating them by $u_0$, so they are no longer emanating from $0$, but from $u_0$. Moreover, taking into account the extra-piece appearing in (\ref{BM kernel}) for $\t_k<\t_\ell$, one computes
  \bea
\lefteqn{ \left. \frac{{dy}}{\sqrt{\pi (t_\ell-t_k)}}
            e^{-\frac{(x-y)^2}{t_\ell-t_k}}
            e^{\frac{x^2}{1-t_k}-\frac{y^2}{1-t_\ell}}
          \right|_{\footnotesize\begin{array}{l}
    t_i= t_0+ (c_0\mu)^2\frac{2\tau_i}
{n^{1/2}}
\\
x= c_0\left( {z_0}  n^{1/2}+ A
\tau_k+\mu\frac{\xi}{n^{1/4}}\right)\\
 y= c_0\left( {z_0}  n^{1/2}+ A
\tau_\ell+\mu\frac{\eta}{n^{1/4}}\right)\end{array}}  }\no\\\no\\
 &=&
    e^{
  {  u_0\mu n^{1/4} (\xi-\eta)} }e^{ \frac{1}{2}\sqrt{n}u_0^2\mu^2(\tau_k-\tau_\ell) }
  e^{ \frac{1}{2}t_0u_0^2\mu^4(\tau_k^2-\tau_\ell^2)}
  \no\\ \no\\
  &&\frac{{d\eta}}{\sqrt{2\pi (\tau_\ell-\tau_k)}}
e^{-\frac{(\xi-\eta)^2}{2(\tau_\ell-\tau_k)}}\left( 1+O(n^{-1/4})\right)
 . \label{rescaling}\eea
  In other terms, setting
 $$D(\xi,\t):=\diag\left(\ldots,e^{
  { -u_0\mu \xi} n^{1/4}}e^{-\frac{1}{2}\sqrt{n}\tau_k u_0^2\mu^2}
  e^{-\frac{1}{2}t_0u_0^2\mu^4\tau_k^2},\ldots\right)_{1\leq k\leq m}
   ,$$
one has the following limit, upon using the $(t_i,x,y)$-rescaling in (\ref{rescaling}) and upon conjugation of the matrix kernel,
 $$\lim_{n\rg \iy}
  D(\xi,\t)(K_{t_k,t_\ell}(x,y))_{1\leq k,\ell\leq m}D(\eta,\t)^{-1}
  =  \left(K^{\cal P}_{\t_k,\t_\ell}(\xi,\eta)\right)_{1\leq k,\ell\leq m},$$
 %
%
%
which leads to the desired kernel (\ref{Pkernel}), upon replacing the integration variables $\om_u\rg U,~\om_v\rg V$.

Since the above argument is obviously formal, one needs to make a rigorous steepest descent analysis on the conjugated kernel above. In the following section, steepest descent contours will be found, depending on whether $q<1,~q=1~q>1$, or what is the same $p>\frac{1}{2},~p=\frac{1}{2},~p<\frac{1}{2}$. Notice the duality $q<1 \leftrightarrow q>1$.


\section{Steepest descent analysis}

In this section, we will deform the contours for both the $u$ and $v$ integration into steepest decent contours. They will be as depicted in the picture below; all lines are at an angle of $0, \pi/4$ or $\pi/2$ with the horizontal line.   

After having set $U:= \frac{c_0u\sqrt{n}}{t_0}$ and $V:= \frac{c_0v\sqrt{n}}{t_0}$ in the exponential appearing in the kernel, one was led to a function
    $$
    F(u) =\frac{u^2}{2}-uz_0+p\log(u-\a)
+(1\!-\!p)\log(u- \b),$$
where
(upon using (\ref{double root}), (\ref {alpha-beta}) and (\ref{u0})),
$$ z_0=\frac{x_0}{c_0}=\b+\frac{2q-1}{r}=\a+\frac{q-2}{r}, ~\a=\frac{at_0}{c_0},~\b=\frac{bt_0}{c_0}, ~u_0=\b+\frac{q}{r}
.$$
 Then one checks (with $p=(1+q^3)^{-1}$ and $r=\sqrt{q^2-q+1}$),
 \be
  z_0-u_0=\frac{q-1}{r},~\qquad u_0-\a=-\frac{1}{r},~\qquad u_0-\b=\frac{q}{r}.
\label{Id}\ee

  \hspace{-2.46cm}
\begin{picture}(0,151.6) 
\put(220,70){\line(1, 1){50 }}
\put(220,70){\line(-1, -1){50 }}
\put(220,19){\line(0, 1){120}}
%
 \put(220,70){\line(1,-1){50}}
  \put(220,70){\line(-1,1){50}}

\put(170,3){\makebox(0,0)[b]{$v$}}
\put(220,3){\makebox(0,0)[b]{$u$}}
\put(270,3){\makebox(0,0)[b]{$v$}}

 \put(170,20){\vector(1, 1){30}}
 \put(220,90){\vector(0, 1){30}}

  \put(210,80){\vector(-1, 1){30}}
  \put(215,75){\vector( 1, -1){30}}
  \put(265,115){\vector( -1, -1){30}}

 \put(105,75){\makebox(0,0)[b]{$q=1$}}

\end{picture}



 \vspace*{4cm}

  \vspace*{ 1cm}\hspace{-2.46cm}
\begin{picture}(12,1.6) 
\put(220,70){\line(1, 1){50 }}
\put(220,70){\line(-1, -1){50 }}
\put(220,19){\line(0, 1){120}}
%
\put(270,120){\line(1, 0){120}}
\put(270,19){\line(1, 0){120}}
 \put(220,70){\line(1,-1){50}}
  \put(220,70){\line(-1,1){50}}

   \put(170,20){\vector(1, 1){30}}
   \put(220,90){\vector(0, 1){30}}
    \put(210,80){\vector(-1, 1){30}}
    \put(215,75){\vector( 1, -1){30}}
    \put(265,115){\vector( -1, -1){30}}

   \put(100,80){\makebox(0,0)[b]{$q>1$}}
\end{picture}


\vspace*{3cm}

  \vspace*{ 1cm}\hspace{-2.46cm}
\begin{picture}(12,1.6) 
 \put(50,120){\line(1, 0){120}}
 \put(50,19){\line(1, 0){120}}
\put(220,70){\line(1, 1){50 }}
\put(220,70){\line(-1, -1){50 }}
\put(220,19){\line(0, 1){120}}
%
 \put(220,70){\line(1,-1){50}}
  \put(220,70){\line(-1,1){50}}
 \put(100,80){\makebox(0,0)[b]{$q<1$}}

  \put(170,20){\vector(1, 1){30}}
  \put(220,90){\vector(0, 1){30}}
  \put(210,80){\vector(-1, 1){30}}
  \put(215,75){\vector( 1, -1){30}}
   \put(265,115){\vector( -1, -1){30}}

 \end{picture}



{\bf (i)}  First we derive the steepest descent contour for the $q$-independent {\bf $u$-integration}, which is the vertical line through $u_0$. Indeed, one checks that for
$$ {\Re e F(u_0+iy)} =\frac{1}{2} (u_0^2-y^2)-u_0z_0+\frac{p}{2}\log ((u_0-\a)^2+y^2)
+(1-p)\log  ((u_0-\b)^2+y^2)
$$
and using (\ref{Id}), the derivative equals
\bean
\frac{\pl}{\pl y} {\Re e F(u_0+iy)}
 &=& -y^3\frac{y^2+\frac{q}{r^2}}{(y^2+\frac{1}{r^2})(y^2+\frac{q^2}{r^2})}
   =-y^3 \left\{ \begin{array}{l}\mbox{a positive}\\ \mbox{ function}\end{array}\right\}
,\eean
showing $ {\Re e F(u_0+iy)} $ has a maximum at $y=0$, which takes care of the $u$-contour.

{\bf (ii)} The next point is to deal with the {\bf $v$-integration}.
Returning to the non-intersecting
Brownian motion kernel (\ref{BM-kernel}), the $F$-function goes with a negative sign on the $v$-contour and thus (with $\vr=\pm 1$)
%
%
 %
\bean
  {-\Re e F(u_0+x(\vr+i))}
 &=&-\frac{1}{2}((u_0+x\vr)^2-x^2)+(u_0+x\vr)z_0\\
 &&-\frac{p}{2}\log((u_0-\a+x\vr)^2+x^2)
    \\
    &&  -\frac{1-p}{2}\log((u_0-\b+x\vr)^2+x^2)
  \eean


 \vspace*{ 6cm}\hspace{-2.46cm}
\begin{picture}(12,1.6) 
\put(50,120){\line(1, 0){120}}
\put(50,19){\line(1, 0){120}}
\put(220,70){\line(1, 1){50 }}
\put(220,70){\line(-1, -1){50 }}
%
\put(270,120){\line(1, 0){120}}
\put(270,19){\line(1, 0){120}}
 \put(220,70){\line(1,-1){50}}
  \put(220,70){\line(-1,1){50}}
\end{picture}

\vspace*{-3.4cm}

$$ \hspace*{-1cm}
 u_0$$

 \vspace*{-3.5cm}

$$\hspace*{6cm}u_0+s(1+i)+x
$$



\vspace*{3.3cm}

$$\hspace*{6cm}u_0+s( 1-i)+x$$

\vspace*{-1.3cm}

$$\hspace*{-9cm}u_0-s(1+i)-x$$

\vspace*{-5.7cm}

$$\hspace*{-9cm} u_0-s( 1-i)-x$$

\vspace*{5cm}

$$0<q<1 \hspace*{6cm} 1<q$$

\vspace*{.5cm}

$$s= \frac{q}{r|q-1|}
$$

\vspace*{1.5cm}

\noindent
 Using (\ref{Id}),
  one checks for $\vr=\pm 1$
\bean
\lefteqn{\hspace*{-4cm}- \left(  \bigl(  u_0- \a+x\vr \bigr) ^{2}+{x}^{2} \right)
 \left(  \bigl( {  u_0}-\b+x\vr  \bigr) ^{2}+{x}^{2} \right)
\frac{\pl}{\pl x} \Re e F(u_0+x(\vr +i))}\\
&=&
-\frac{4 {x}^{3}}{r} \left( { { -\left( q-1 \right) \vr x}}+{\frac {q}{r}} \right)
\\
&=&  -x^3\times \left\{ \begin{array}{l}\mbox{a positive}\\ \mbox{ function}\end{array}\right\}
\eean
for
\be
\begin{array}{ll}
   \left\{\begin{array}{ll}-\iy< x< \frac{q}{r(q-1)} ,&~\vr=1\\
     - \frac{q}{r(q-1)}<x<\iy,&~\vr=-1
   \end{array}\right\}&~~\mbox{if} ~q>1,
 \\  \\
    \left\{\begin{array}{ll}-\iy< x< \iy,&~\vr=1\\
     -\iy<x<\iy,&~\vr=-1
   \end{array}\right\}&~~\mbox{if} ~q=1,
 \\  \\
   \left\{\begin{array}{ll} -\frac{q}{r(1-q)}<x<\iy ,&~\vr=1\\
   -\iy< x < \frac{q}{r( 1-q)},&~\vr=-1
   \end{array}\right\}&~~\mbox{if} ~q<1,
 \end{array}
 \label{segments}\ee
 Therefore the function  $- \Re e F(u_0+x(\vr +i))$, restricted to the segments specified by (\ref{segments}),  has its maximum at $u_0$. One then completes those segments by horizontal lines starting from the end of those segments as in the figure above. Along those horizontal half lines, one must check that the maximum is attained at the points $u_0\pm s(1-i)$ and $u_0\pm s(1+i)$.
To carry out this computation,
the four horizontal segments can readily be represented by
  $$u_0+\dt s(\vr+i)+\dt \vr x$$
   with
  $$\dt=-\vr=1\qquad\qquad\qquad\qquad \dt=\vr=1$$

  \begin{picture}(1250,0)  \hspace*{6.6cm}
\put(0,0){\line(1, 1){25}}
\put(0,0){\line(-1, -1){25}}
\put(1,1){\line(-1, 1){25}}
\put(1,1){\line( 1, -1){25}}
\end{picture}

  $$ -\dt= \vr=1\qquad\qquad\qquad\qquad -\dt=-\vr=1$$
  in the four corresponding regions of the figure above.


 To deal with the horizontal segment, setting $s=\frac{q}{r|q-1|}$, one finds 
 \bean
  {-\Re e F(u_0+\delta s(\vr+i)+\delta\vr x)}
 &=&-\frac{1}{2}((u_0+\epsilon\delta(s+x))^2-s^2)+(u_0+\epsilon\delta(s+x))z_0\\
 &&-\frac{p}{2}\log\left((u_0-\a+\epsilon\delta(s+x))^2+s^2\right)
    \\
    &&  -\frac{1-p}{2}\log\left((u_0-\b+\epsilon\delta(s+x))^2+s^2\right)
 . \eean
  %
  One computes,
  \bean
  &&-
  \left(  \bigl( {  u_0}+(s+x)\epsilon  \delta-{ \a}
 \bigr) ^{2}+s^{2} \right)  \left(  \bigl( {  u_0}+(s+x)\epsilon
\delta -{  \b} \bigr) ^{2}+{s}^{2} \right)
\\
&&\hspace*{5cm}\frac{\pl}{\pl x}\Re e ~F(u_0+\delta s(\vr+i)+\delta\vr x)
 \Bigr|_{s=  \frac{\vr \dt q}{r(q-1)}}
\eean
  \bea
  = \frac{q^4}{r^5(q-1)^5}\left.\left(\begin{array}{l}
   q z^5 +\vr \dt (q^2+3q+1)z^4+3   (q+1)^2z^3 + \vr \dt  (q+3)z^2 \\
   +4    ( {q}^{2}+q+1  ) z +\vr \dt  ( 1+{q}^{2} )
   \end{array}\right)\right|_{z=\frac{r(q-1)}{q}x}\label{3.3}\eea
 Since $x>0$ is increasing as one moves away from $u_0$ along the horizontal lines; one has simultaneously,
 $$
 \begin{array}{lllll} s=\frac{q}{r(q-1)},& \vr \dt=+1,& q-1>0,&z>0\\
 s=-\frac{q}{r(q-1)},& \vr \dt=-1,& q-1<0,&z<0\end{array}
  $$ and thus, since $q>0$, the right hand side of (\ref{3.3}) is $>0$. Thus the above derivative is negative on the four lines and so, when one moves away from $u_0$ along the horizontal paths in all four directions (as in the picture above) the function $-\Re e~F(z)$ goes down so that combining both calculations, the maximum will be attained at $u_0$.

In order to show that in the limit the Pearcey kernel is obtained, one picks $\tau_i$'s, and $\xi, \eta$ in a compact set of $\BR$, one integrates the $u$ and $v$ variables along the contour in a neighborhood of radius $n^{-1/4} n^{\frac{1}{20}-\vr}$. Then $\left|
 \frac{\om_u}{\mu}\right|$,~$\left|
 \frac{\om_u}{\mu}\right|\leq \dt\leq n^{\frac{1}{20}-\vr}
 $, as $n\rg \iy$. In this range lemma \ref{tech} will apply, whereas outside this neighborhood, the rest of the contour makes no contribution, because of the steepest descent estimates.
So, one needs the following estimate:
\begin{lemma}\label{tech}

Given the function $F(u)$ as in (\ref{F-function}), one has the following estimate,
$$n\left |
  F(u_0+\frac{\dt}{n^{1/4}})-F(u_0)-F^{(iv)}(u_0)\frac{1}{4!} \left(\frac{\dt}{n^{1/4}}\right)^4
 \right|
 \leq \frac{64\dt^5}{5n^{1/4}} \left(q+\frac{1}{q}\right)^5
  .$$
\end{lemma}

\proof By Taylor's Theorem and using $
  F'(u_0)=F''(u_0)=F'''(u_0)=0  \mbox{  and  }
  \frac{1}{4!}F^{(iv)}(u_0)=   -\frac{r^2}{4q}
   $, as in (\ref{derivFnul}), one has
 $$
  \left |
  F(u_0+\frac{\dt}{n^{1/4}})-F(u_0)-\frac{1}{4!} \left(\frac{\dt}{n^{1/4}}\right)^4F^{(iv)}(u_0)
 \right| \leq  \left(\frac{\dt}{ n^{1/4}}\right)^5 \max_{ |u-u_0|\leq \frac{\dt}{n^{1/4}}}
  \frac{\left| F^{(v)}(u) \right|}{5!}
 $$
 From the explicit expression (\ref{F-function}) for $F$, from the fact that $u_0-\b=q/r$ and $a-u_0=1/r$ and that $\beta<u_0<\a$, one deduces\footnote{Remember $r=\sqrt{q^2-q+1}$}
  \bean
  \max_{|u-u_0|\leq \frac{\dt}{ n^{1/4}}}\left|\frac{ F^{(v)}(u)}{5!}\right|
   &=& \sup_{|u-u_0|\leq \frac{\dt}{ n^{1/4}}}
   \frac{1}{5} \left| \frac{1-p}{(u-\b)^5}+
   \frac{p}{(u-\a)^5}\right|   \\
   &\leq&
    \frac{2}{5   \min (| \a-u_0- \frac{\dt}{ n^{1/4}} |,| u_0-\b- \frac{\dt}{ n^{1/4}} |)^5}
   \\
   &\leq &
    \frac{2}{5}
      \left( \frac{1}{r}\min(1,q)- \frac{\dt}{ n^{1/4}}\right)^{-5}
      \\
      &\leq&
       \frac{64}{5}\left(\frac{r }{(\min(1,q)) }\right)^5
        \left(\begin{array}{l}
         \mbox{  by picking $n$ large enough}\\
         \mbox{ such that } \frac{\dt}{ n^{1/4}}\leq  \frac{1}{2r}\min(1,q)
       \\ \mbox{since $\dt\leq n^{\frac{1}{20}-\vr}$}  \end{array}\right)\\
         &\leq &
          \frac{64}{5}\left(q+\frac{1}{q}\right)^5
     \eean
ending the proof of Lemma \ref{tech}.\qed


\section{Proof of Theorem \ref{Theo:2}}





In this section, we denote by $a$ and $b$ the Brownian
motions target points, where we put $b=0$. We denote the
old time 
in the Brownian motion formula by $\bar t$ and the new rescaled
time and space in the Brownian motion formula (Theorem
\ref{Theo 1}) by $\bar \tau$ and $\bar \eta$. Set $   E^c=\cup_{i=1}^r(y_{2i-1},y_{2i})\subset \BR$. Let $x$ be the spatial variable for the
matrix integral and $\a$ the variable appearing in the
diagonal matrix, the other one being $=0$.

The reader is reminded of the different players in the argument below, in accordance with (\ref{Brownian transition}),
 \bean  
\lefteqn{\BP_{\mbox{\tiny Br}}(\bar t, y, a\sqrt{n})}\\
&&:=\BP^{(0,a\sqrt{n})}_{n} \left(\begin{tabular}{l|l}
& all $x_j(0) =0$\\
all $x_j(\bar t)\in E$ for $1\leq j\leq n$ &
 $n_1$ paths end up at $a\sqrt{n}$ at $\bar t=1$\\
  & 
$n_2$ paths end up at $0$ at $\bar t=1$
\end{tabular}\right)
\eean
and, setting $\tilde b_1=\a,~\tilde b_2=0$,
 $$
\BP_n(\a,x_i):=\frac{1}{Z_n}
 \int_{\tilde E^n}
 \Delta_n(x_1,...,x_n)
%
 \prod^2_{\ell=1} \Delta_{n_{\ell}}
(x^{(\ell)})\prod_{j=1}^{n_{\ell} }
e^{-\frac{1}{2}x_j^{(\ell)^2}+\tilde
b_{\ell}x_j^{(\ell)}
 }dx_j^{(\ell)}
$$
 with $\BP_n(\a,x_i)$ satisfying the PDE\footnote{ Given $ \tilde E^c=\cup_{i=1}^r(x_{2i-1},x_{2i})\subset \BR$, the prime in the formula below denotes $':= \sum \frac{\pl}{\pl x_i}$.}, as in (\ref{7.4.11}) and (\ref{eqts}),
  $$
\det\left(\begin{array}{ccc}
 F_1 & F_2 & 0\\
 F_1'& F_2'&
 F_1F_2\left(\frac{H_1}{F_1}+\frac{H_2}{F_2}\right)\\
 F_1''& F_2''&
 F_1F_2\left(\frac{H_1'}{F_1}+\frac{H_2'}{F_2}\right)
\end{array}\right)=0
 ,$$
  with $F_i$ and $H_i$ given in (\ref{F,G,H}).
  From (\ref{Brownian transition}), one has the following relationship
  $$
 \BP_{\tiny Br}\left(\bar t, y, a\sqrt{n}\right)
  =\BP_{n}\left(a\sqrt{n}\sqrt{\frac{2\bar t}{1-\bar t}},
  y_i\sqrt{\frac{2}{\bar t(1-\bar t)}}~\right)=\BP_n(\a,x_i)
 $$
 and also from Theorem \ref{Theo 1},
 $$
  \left. \BP_{\tiny Br}(\bar t, y, a\sqrt{n})\right|_{\bar t=t_0+ (c_0\mu)^2\frac{2\bar\tau}
{n^{1/2}} ,~~ y= {x_0}  n^{1/2}+c_0 A
\bar\tau+c_0\mu\frac{\bar \eta}{n^{1/4}}   }
 $$
\be =\BP^{\cal P}\left ({\cal P}(\bar \tau)\cap \cup_{i-1}^r(\bar \eta_{2i-1},\bar \eta_{2i}) =\emptyset\right)+O(n^{-1/4})
\label{approx}\ee
Setting $b=0$, the formulae (\ref{notation}), (\ref{x0,t0}) and (\ref{constants}) in Theorem \ref{Theo 1} simplify: 
$$
 \sqrt{\frac{2t_0}{1-t_0}}=\frac{q+1}{ ar},~~
  x_0=
  \frac{2q-1}{q+1}at_0
  ~,~~c_0=\sqrt{\frac{t_0(1-t_0)}{2}}= \frac{at_0 r}{q+1}$$
\be
A =
\frac{q -(2q-1)
t_0}{\sqrt{q} },~~~\mu=\frac{\sqrt{r}}{q^{1/4}}.\label{simple}\ee
For convenience, we set
 $$  \tau:=  \frac{\bar \tau}{\sqrt{q}}, ~~
\eta :=\frac{\bar \eta}{q^{1/4}},~~z:=\left(\frac{q^2-q+1}{n}\right)^{1/4}.$$
One finds, using $z^2=r/\sqrt{n}$ and formulae (\ref{simple}), 
\bean
 \bar t= t_0+(c_0\mu)^2\frac{2\bar \tau}{\sqrt{n}}
 &=&
   t_0\left(1+(1-t_0)\frac{\bar \tau}{\sqrt{q}}
   \left(\frac{q^2-q+1}{n}\right)^{1/2}\right)
    \\
   &=&
     t_0\left(1+(1-t_0) \tau z^2\right)\label{t}\eean
 Moreover, one computes
\bean
 \a =a\sqrt{n}\sqrt{\frac{2\bar t}{1-\bar t}}
 &=&\sqrt{n}a\sqrt{\frac{2t_0}{1-t_0}}
  \sqrt{\frac{1+(1-t_0)  \tau z^2}{1-t_0  \tau z^2}}
  %
   %
     \\
      &=&
    \frac{q+1}{z^2}
     \sqrt{
     \frac{1+(1-t_0)  \tau  z^2
     }{1-t_0  \tau z^2}
     }.
\eean
 We also have, using the formulae (\ref{simple}) above,
 \bean
 x &=&y\sqrt{\frac{2}{\bar t(1-\bar t)}}\\
 &=&
 \sqrt{\frac{2}{\bar t(1-\bar t)}}
c_0 \left(\frac{x_0}{c_0}\sqrt{n}+ A\bar \tau+ \mu\frac{\bar \eta}{n^{1/4}}\right)
 \\&=&
 \sqrt{\frac{2}{\bar t(1-\bar t)}}
  c_0\left(
  (2q-1)\frac{\sqrt{n}}{r}+(q-(2q-1)t_0)\frac{\bar \tau}{q^{1/2}}
  +\frac{\bar \eta}{q^{1/4}}
  \frac{\sqrt{r}}{n^{1/4}}\right)
  %
 \\
 &=&\frac{1}{\sqrt{(1+(1-t_0) { \tau}  z^2
     )}\sqrt{( {1-t_0 \tau z^2}
     )}}\left(
      \frac{2q-1}{ z^2}+(q-t_0(2q-1))
      \tau+   \eta z\right)\eean
One also checks
 $$n_1=np=\frac{n}{ q^3\!+\!1 }=\frac{n}{(q\!+\!1)(q^2-q+1)}=
  {\frac { {z}^{-4}}{ \left( q+1 \right) }},~~n_2=n(1-p)=
 {\frac {{q}^{3}{z}^{-4}}{ \left( q+1 \right)  }}
   .$$
%
Consider the map $T_z:(\tau,\eta_i)\mapsto (\a, x_i)$,
 \bean
 (\a,x_i)&:=&T_z(\t,\e_i)
  \\
  &=& \left( {\frac {( q+1) \sqrt{
   1+  ( 1-{ t_0})
   \t{z}^{2}  }    }{{z}^{2}
 \sqrt {\left( 1-{ t_0} \t{z}^{2}
\right) }}},    {\frac {\left( {\frac {2q-1}{{z}^{2}}}+
\left( q-  ( 2q-1  ) {   t_0} \right) \t+\e_i z
\right)}{\sqrt { \left( 1+  ( 1-{   t_0} ) \t{z}^{2}
\right)  \left( 1-{ t_0}\t{z}^{2} \right) }}} \right)
\eean
and its inverse $T_z^{-1}: (\a, x_i)\mapsto (\tau,\eta_i)$,
\bean
 (\t,\e_i)&=&T^{-1}_z(\a,x_i)\\
  &=&\left(  \frac{{\a}^{2}{z}^{4}-  ( q+1  )^2
}{z^2({  t_0} {z}^{4}{\a}^{2}+ ( q+1) ^{2}
 ( 1\!- \!t_0
  ))},~  \frac{\a (  x_i q-\a q+x_i )
 {z}^{4}-(q\!-\!1)(q\!+\!1)^2
}{z^3({  t_0} {z}^{4}{\a}^{2}+ \left( q+1 \right) ^{2}
\left( 1-{  t_0}
 \right))}
\right)\eean
Then setting
\bean \log \BP_n(\a,x_i)&=&F(\t,\e_i)=F(T_z^{-1} (\a,x_i)),
 \eean
setting $B=\cup_{i-1}^r(x_{2i-1},x_{2i})\subset \BR$, taking the
derivatives $\pl_x:=\sum \frac{\pl}{\pl x_i},~~\vr_x:=\sum x_i\frac{\pl}{\pl x_i}$ and then taking a
series in $z$, the functions $F_i$ and $H_i$ in (\ref{F,G,H}) have the
following form, 
 \bean
 F_1&=& 
   -\frac{\pl}{\pl \a}\pl_{_x}\log \BP_n+n_1
   =\frac{1}{q+1}\left(z^{-4}+z^{-2}{ F''}
 -2z^{-1}\frac{\pl F'}{\pl \t} +O(1)
\right)\\
 F_2&=& (\frac{\pl}{\pl \a}+\pl_{_x})
  \pl_{_x}\log \BP_n+n_2=
   \frac{1}{q+1}\left(q^3z^{-4}+qz^{-2}  F''
 +2z^{-1}\frac{\pl F'}{\pl \t} +O(1)
\right)
  \\
  H_1&=&\!\left\{H^{(1)}_1,F_1\right\}_{\pl_{_x}}\!\!+
  \left\{H^{(2)}_1,F_1\right\}_{\a}
  =\frac{2}{(q\!+\!1)^4} \left(
   \begin{array}{l}
    q(2q^2+3q+2)F'''z^{-9}
    \\
    - (4{q}^{3}-1+6{q}^{2}+3q)  \frac{\pl F''}{\pl \t}z^{-8}
 \\+O(z^{-7})
 \end{array}
 \!\!\! \right)
\\
H_2&=&\!\left\{H^{(1)}_2,F_2\right\}_{\pl_{_x}}\!\!-
\left\{H^{(2)}_2,F_2\right\}_{\a}
=\frac{2q^3}{(q\!+\!1)^4}\left(\begin{array}{l} - {q}
 ( 2{q}^{2}+3q+2 )F'''  {z}^{-9}\\+
   \left( 4{q}^{3}\!-\!1+\!6{q}^{2}\!+\!3q
 \right) \frac{\pl F''}{\pl \t}{z}^{-8}\\
 +O \left( {z}^{-7} \right)
\end{array}\!\!\!\right)\label{F's}\eean
with
 \bean
H_1^{(2)}&=&
  (1-\vr_{_x}+\a\frac{\pl}{  \pl
\a}+
 2\a\pl_{_x})
 (-\frac{\pl}{  \pl \a})\log \BP_n\\
H_2^{(2)}&=&
 (1-\vr_{_x}+\a\frac{\pl}{  \pl \a})
 (\frac{\pl}{  \pl \a}+\pl_{_x})\log \BP_n%
 \\
 H_1^{(1)}&=&
  -2\frac{\pl}{\pl \a}\log \BP_n-2n_1(\a+\frac{n_2}{\a})\\
  H_2^{(1)}&=&
    (1+\vr_{_x} -\a\frac{\pl}{\pl \a})
    \frac{\pl}{\pl \a} \log \BP_n
    +2\frac{n_1n_2}{\a}
  \eean
Then one computes (\ref{eqts}), setting ${}'={\pl_x}$ and setting
$\pl_E=\sum\frac{\pl}{\pl y_i} $ for $E=(y_1,y_2)$,
$$
\left.\det\left(\begin{array}{ccc}
 F_1 & F_2 & 0\\
 F_1'& F_2'&
 F_1F_2\left(\frac{H_1}{F_1}+\frac{H_2}{F_2}\right)\\
 F_1''& F_2''&
 F_1F_2\left(\frac{H_1'}{F_1}+\frac{H_2'}{F_2}\right)
\end{array}\right)\right|_{(\a,x_i)=T_z(\t,\e_i),
 ~\mbox{with}~\t=\frac{t}{\sqrt{q}},~
  \e_i=\frac{y_i}{{q}^{1/4}}}
~~~~~~~~~~~~~~~~~~~~~~~~~~~~~$$
\bean
&=&-2q^{6+1/2}\frac{q-1}{(q+1)^5} \Bigl\{\pl_{_E}^3 \log
\BP^{\cal P}, \BX   \Bigr\}_{_{\pl_E}}z^{-18}
\\
&&+\frac{1}{8}\left(\left\{\frac{\pl}{\pl t}\pl_{_E}^2
\log
\BP^{\cal P},\BX\right\}_{_{\pl_E}}+O(q-1)\right)z^{-17}+O(z^{-16}),\eean
%
%
where, setting  ${\mathbb Q}(t,E):= \log\BP^{\cal P}\left({\cal P}(t ) \cap E =
\emptyset 
 \right)$,
$$
 \BX:=
  8\frac{\pl^3 {\mathbb Q}}{\pl t^3}
  +\left(\vr_{_{\! E}}-2t\frac{\pl}{\pl t}-2\right)\pl_{_{\! E}}^2{\mathbb Q}
 -4\left\{\pl_{_E}^2 {\mathbb Q},\pl_{_E} \frac{\pl {\mathbb Q}}{\pl t} \right\},$$
where we made use of (\ref{approx}), which states that $\log \BP_n(\a, x_i)={\mathbb Q}(t,E)+O(z)$.

 \noindent$\bullet$ For $q\neq 1$, the function $\log \BP^{\cal P}$, which is independent of $q$ by the universality result,
  satisfies
 the differential equation, given by the leading term
 $z^{-18}$,
 \be \left\{\pl_E^3 \log \BP^{\cal P},\BX\right\}=0.\label{eqt1}\ee

  \noindent$\bullet$ For $q=1$, the $z^{-18}$-term vanishes and thus $\log
 \BP$ satisfies another equation, namely the one appearing in the
 $z^{-17}$-term,
  \be\left\{\pl_E^2\frac{\pl  }{\pl t}\log \BP^{\cal P},\BX\right\}=0
  .\label{eqt2}\ee
  This means that $\log \BP^{\cal P}$ satisfies the
  two equations (\ref{eqt1}) and (\ref{eqt2}).

Thus
 for $E=(x,y)\subset \BR$, setting $u_{\pm}=\frac{1}{2}
 (y\pm x)$, and
 $\log \BP^{\cal P}=H(t;\frac{1}{2}(y+x),\frac{1}{2}(y-x))$,
  \bean
   \pl_{_E}\log \BP^{\cal P}&=&(\frac{\pl}{\pl x}+\frac{\pl}{\pl
  y})\log \BP^{\cal P}
   =\frac{\pl}{\pl u_+}H\\
    \vr_{_E}\log \BP^{\cal P}&=&(x\frac{\pl}{\pl x}+y\frac{\pl}{\pl
    y})\log \BP^{\cal P}
   =(u_+\frac{\pl}{\pl u_+}+u_-\frac{\pl}{\pl u_-})H
    .\eean
 Since the wronskian of two functions equals the derivative
 of the ratio, modulo a non-zero multiplicative term, one
 concludes from equation (\ref{eqt1}) that
 $
 \BX=c(t,u_-) \frac{\pl^3 H}{ \pl u_+^3}
 ,$
  with $c(t,u_-)$ a function depending on all variables except
  $u_+$;
   putting this equation in equation (\ref{eqt1}),
 one finds
  $$
    c(t,u_-)\left\{\frac{\pl^3 H}{\pl u_+^3},
   \frac{\pl^3 H}{\pl t \pl u_+^2}\right\}_{u_{+}}=
    c(t,u_-)\left\{\pl^3_E\log \BP^{\PR} , \frac{\pl}{\pl t}\pl^2_E\log \BP^{\PR}\right\}_{\pl_E}
=0,
   $$
   implying $c(t,u_-)=0$ for all $t, u_-$ and thus
   $H(t;u_+,u_-)=\log \BP^{\PR}$ satisfies the equation $\BX=0$, provided the Wronskian $\left\{\pl^3_E\log \BP^{\PR} , \frac{\pl}{\pl t}\pl^2_E\log \BP^{\PR}\right\}_{\pl_E}\neq 0$. This will be shown in the next section, using functional theoretical arguments. This ends the proof of Theorem \ref{Theo:2}, except for the ``final condition", which will be shown in \cite{AvMCafasso}. \qed


\section{An estimate for the Wronskian}


\begin{proposition}\label{Prop-A1}
The Wronskian

$$
\left\{ \frac{\pl}{\pl t}\pl^2_E\log \BP^{\PR}, \pl^3_E\log \BP^{\PR} \right\}_{\pl_E}
$$
is a non-zero function.
\end{proposition}

In order to prove this proposition, we need a number of lemmas; the proofs will be functional-theoretical and rely on the techniques and on some of the formulae in \cite{TW-Airy}. Given the Pearcey kernel $K^{\PR}$, one defines, for a given set $E=\cup_{k=1}^r [a_{2k-1},a_{2k}]\subset \BR$, the following\footnote{Given a kernel, viewed as an operator, the equality $\doteq$ refers to the corresponding kernel}:
\be
K_E^{\PR} :=K^{\PR}\raisebox{1mm}{$\chi$}{}_{E},\qquad I+R:=(I-K_E^{\PR})^{-1}\doteq :\rho(x,y)
\label{6.notation} \ee
and thus one has
 \be
 (I-K_E^{\PR})^{-1}K_E^{\PR}=K_E^{\PR}+(E_E^{\PR})^2+... =R.
\label{R}\ee
Also, from (\ref{Pkernel}) and the differential equations (\ref{ode}) for $p(x)$ and $q(y)$, it follows that
\be
\left(\frac{\pl}{\pl x}+ \frac{\pl}{\pl y}\right)K^{\PR}(x,y)=p(x)q(y).
\ee
Note that for a general kernel $L$, one has

\be
[D,L]\doteq\left(\frac{\pl}{\pl x}+ \frac{\pl}{\pl y}\right)L(x,y).
\ee
Define the functions
\be
\hat p:=(I-K_E^{\PR})^{-1}p,\quad \hat q:=(I-K_E^{\PR^{\top}})^{-1}q
\label{p hat}\ee
and\footnote{$\la f,g\ra:= \int_{\BR} f(x)g(x)dx$.}
\be
u:=\int_E(I-K_E^{\PR})^{-1}p(\alpha)q(\alpha)d\alpha=\la\hat p(\alpha),q(\alpha)\raisebox{1mm}{$\chi$}{}_{E}(\alpha)\ra .
\label{A.41}\ee

\begin{lemma} \label{Lemma 6.2}For a disjoint union $E=\bigcup^r_{k=1}[a_{2k-1},a_{2k}]$, one has the following identity\footnote{Remembering $\pl_E=\sum_1^r \frac{\pl}{\pl a_j}$.}:
$$
\pl^2_E\log \BP^{\PR}=\pl^2_E\log\det(I-K_E^{\PR})=\pl_Eu=\sum_k (-1)^k\hat p(a_k)\hat q(a_k).
$$
\end{lemma}

\proof At first notice that
\bea
[D,K_E^{\PR}]=[D,K^{\PR}\raisebox{1mm}{$\chi$}{}_{E}]&\doteq&
  \left(\frac{\pl}{\pl x}+ \frac{\pl}{\pl y}\right)(K^{\PR}\raisebox{1mm}{$\chi$}{}_{E})
  \label{A.5}\\
\nonumber\\
&=&p(x)q(y)\raisebox{1mm}{$\chi$}{}_{E}(y)-\sum_k(-1)^kK^{\PR}(x,a_k)\dt(y-a_k),
\nonumber\eea
from which one deduces, using notations (\ref{6.notation}) and (\ref{R}),
\bea
\left(\frac{\pl}{\pl x}+ \frac{\pl}{\pl y}\right)R(x,y)&\doteq&[D,(I-K_E^{\PR})^{-1}]\nonumber\\
\nonumber\\
&\doteq&(I-K_E^{\PR})^{-1}p(x)q(y)\raisebox{1mm}{$\chi$}{}_{E}(y)(I-K_E^{\PR})^{-1}\nonumber\\
\nonumber\\
& &-\sum_k(-1)^k(I-K_E^{\PR})^{-1}K^{\PR}(x,a_k)\dt(y-a_k)(I-K_E^{\PR})^{-1}\nonumber\\
\nonumber\\
&=&\hat p(x)\hat q(y)\raisebox{1mm}{$\chi$}{}_{E}(y)-\sum _k(-1)^kR(x,a_k)\rho(a_k,y).
\label{A.6}\eea
and
\bea
\frac{\pl}{\pl a_k}R(x,y)= \frac{\pl}{\pl a_k}  (I+R)
&=&\frac{\pl}{\pl a_k}(I-K_E^{\PR})^{-1}\nonumber\\
\nonumber\\
&=&(I-K^{\PR}\raisebox{1mm}{$\chi$}{}_{E})^{-1}K^{\PR}\frac{\pl\raisebox{1mm}{$\chi$}{}_{E}}{\pl a_k}(I-K^{\PR}\raisebox{1mm}{$\chi$}{}_{E})^{-1}\nonumber\\
\nonumber\\
&=&R(x,z)(\dt(z-a_k)(-1)^k)\rho(z,y)\nonumber\\
\nonumber\\
&=&(-1)^kR(x,a_k)\rho(a_k,y).
\label{A.7}\eea
Then combining (\ref{A.6})  and (\ref{A.7}), one finds
\bea
\left(\frac{\pl}{\pl x}+ \frac{\pl}{\pl y}+\sum_k \frac{\pl}{\pl a_k}\right)R(x,y)&=&[D,(I-K_E^{\PR})^{-1}]+\sum_k\frac{\pl}{\pl a_k}R(x,y)\nonumber\\
\nonumber\\
&\doteq&\hat p(x)\hat q(y)\raisebox{1mm}{$\chi$}{}_{E}(y),
\eea
and hence, setting $x=y=a_j$, the total derivative becomes
\be
\sum_k\frac{d}{da_k}R(a_j,a_j)\doteq\hat p(a_j)\hat q(a_j).
\label{A.9}\ee
\newcommand{\KEP}{K_E^{\PR}}
  One then computes the derivative of $u$, as defined in (\ref{A.41}), with respect to $a_k$: (of course, the functions $p$ and $q$ do not involve the interval $E$)
\bea
\frac{\pl u}{\pl a_k}&=&\frac{\pl  }{\pl a_k}\Big\la (I-\KEP)^{-1}p,q\raisebox{1mm}{$\chi$}{}_{E}\Big\ra\nonumber\\
\nonumber\\
&=&\Big\la \bigl( \frac{\pl  }{\pl a_k}(I-\KEP)^{-1}\bigr)p,q\raisebox{1mm}{$\chi$}{}_{E} \Big\ra +\Big\la\hat p,q\frac{\pl\raisebox{1mm}{$\chi$}{}_{E} }{\pl a_k}\Big\ra\nonumber\\
\nonumber\\
&=&(-1)^k\la R(x,a_k)\hat p(a_k)
 ,q\raisebox{1mm}{$\chi$}{}_{E}\ra
 +\la\hat p,q\dt(y-a_k)(-1)^k\ra, ~\mbox{  using (\ref{A.7}),}\nonumber\\
\nonumber\\
&=&(-1)^k\left(\hat p(a_k)\la R(x,a_k),q\raisebox{1mm}{$\chi$}{}_{E}\ra +\hat p(a_k)q(a_k)\right)\nonumber\\
\nonumber\\
&=&(-1)^k\hat p(a_k)((I+R^{\top})q(a_k)-q(a_k)+q(a_k))
\nonumber\\
&=&(-1)^k\hat p(a_k)\hat q(a_k).
\label{A.10}\eea
Combining (\ref{A.9}) and (\ref{A.10}) yields
$$
\frac{\pl u}{\pl a_j}=\sum_k\frac{d}{d a_k}(-1)^jR(a_j,a_j)
$$
and then summing with respect to $j$,\footnote{Here one uses identity (1.1) in \cite{TW-Airy},
$$
\frac{d}{da_j}\log\det(I-\KEP)^{-1}=(-1)^{j+1}R(a_j,a_j).
$$}
\bean
\sum_k\frac{\pl u}{\pl a_k}&=&\sum_k\frac{d }{d a_k}\left( \sum_j(-1)^jR(a_j,a_j)\right)\\
&=&-\left(\sum_k\frac{d }{d a_k}\right)^2\log\det(I-\KEP)^{-1}\\
&=&\pl_E^2\log\det(I-\KEP)=\pl_E^2\log \BP^{\PR},
\eean
which, together with (\ref{Pkernel}), establishes Lemma \ref{Lemma 6.2}.\qed

\begin{lemma}\label{Lemma-A2}
Given $E=[x,y]$, the following estimates hold
\bean
 \left( \frac{\pl}{\pl x}+ \frac{\pl}{\pl y}\right)u&=&(y-x)(p(x)q(x))'+{\bf O}(y-x)^2\\
\\
\frac{\pl u  }{\pl t}&=&\frac{1}{2}(y-x)(p''q-pq'')(x)+{\bf O}(y-x)^2.
\eean
\end{lemma}

\proof Using the fact that, for a small interval $E$, the integral $\int_x^y$ has order $x-y$, using $R(\a,y)-R(\a,x)={\bf O}(y-x)$,  one deduces from the formula of Lemma \ref{Lemma 6.2},
(remember the definitions (\ref{p hat}) of $\hat p$ and $\hat q$)
\bean
\pl_Eu&=&\hat p(y)\hat q(y)-\hat p(x)\hat q(x)\\
\\
&=&\left(p(y)+\int_x^yR(y,\a)p(\a)d\a\right)\left(q(y)+\int_x^yR(\a,y)q(\a)d\a\right)\\
\\
& &-\left(p(x)+\int_x^yR(x,\a)p(\a)d\a\right)\left(q(x)+\int_x^yR(\a,x)q(\a)d\a\right)\\
&=&p(y)q(y)-p(x)q(x)+{\bf O}(y-x)^2\\
&=&(y-x)(p(x)q(x))'+{\bf O}(y-x)^2.
\eean
Using the heat equations (\ref{heat}) satisfied by $p$, $q$ and the PDE (\ref{kernel-PDE}) for $K^{\PR}$,
%
one checks
\bean
2\frac{\pl u}{\pl t}&=&2\frac{\pl }{\pl t}\la (I-\KEP)^{-1}p,q\ra \\
\\
&=&\Big\la (I-\KEP)^{-1}2\frac{\pl \KEP}{\pl t}(I-\KEP)^{-1}p,q\Big\ra \\
\\
& &+\Big\la (I-\KEP)^{-1}2\frac{\pl p}{\pl t},q\Big\ra +\Big\la (I-\KEP)^{-1}p,2\frac{\pl q}{\pl t}\Big\ra \\
\\
&=&\Big\la \int_E dy\left(-p'(x)q(y)+p(x)q'(y)\right)\hat p(y),\left((I-K_E^{\PR^{\top}})^{-1}q\right)(x)\Big\ra  \\
\\
& &-\Big\la p''(x),\left((I-K_E^{\PR^{\top}})^{-1}q\right)(x)\Big\ra +  \Big\la(I-\KEP)^{-1}p(x),q''(x)\Big\ra \\
\\
&=& u(x)(-\la p',\hat q\ra+\la q',\hat p\ra ) -\la p'',\hat q\ra +\la\hat p,q''\ra ,\mbox{  using $\la p,\hat q\ra=\la \hat p, q\ra=u$},\\
\\
&=&{\bf O}(y-x)^2-\la p'',q\ra +\la p,q''\ra  
\\
&=&-(y-x)(p''(x)q(x)-p(x)q''(x))+{\bf O}(y-x)^2,
\eean
thus ending the proof of Lemma \ref{Lemma-A2}.\qed

\medskip\noindent{\it Proof of Proposition \ref{Prop-A1}:\/} One computes the following Wronskian and expand for small $y-x$, given $E=(x,y)$. Indeed, from $\pl_E u=\pl^2_E\log \BP^{\PR} $(see Lemma \ref{Lemma 6.2}), the estimates of Lemma \ref{Lemma-A2} and the the differential equations (\ref{ode}) for $p$ and $q$, one computes
\bean
\lefteqn{\left\{ \frac{\pl}{\pl t}\pl^2_E\log \BP^{\PR}, \pl^3_E\log \BP^{\PR}\right\}_{\pl_E}}\\
\\
&=&\left\{\pl_E \frac{\pl u}{\pl t},\pl^2_E u\right\}_{\pl_E}\\
\\
&=&-\frac{(y-x)^2}{2}\left( \{(p''q-pq'')',(pq)''\}_x+{\bf O}(y-x)\right)\\
\\
&=&\frac{(y-x)^2}{2}\left(\begin{array}{c}
(-t(p''q-pq'')+3x(pq)'+2pq)(pq)''+{\bf O}(y-x)\\
+(t(p'q-pq')-2xpq+p''q'-p'q'')(t(pq)'+3(p'q')')
\end{array}\right)
\eean
with the coefficient of $(y-x)^2/2$, for $x=t=0$ being equal to

$$
2pq(pq)''-3(p'q')'(p'q''-p''q')\neq 0
$$
which is nonzero, ending the proof of Proposition \ref{Prop-A1}.\qed

\section{Steepest descent analysis and replica duality}\label{section 7}

In this section, we emphasize the central role played by the spectral curve (or Pastur equation \cite{Pastur}) in the steepest
descent analysis used to prove the universal behavior of the kernel as $N \to \infty$. More precisely, we point out that
the kernel used for the steepest descent analysis takes a very universal form in the different problems of
non-intersecting Brownian motions studied so far. We further give show that the study of the spectral curve associated
to the considered problem gives the universal limiting kernel obtained after a proper rescaling of the variable
around a singularity of the problem: we give a "physical meaning" to the computations performed in the preceding
sections as well as a way to generalize it to more complicated problems.

In a first part, we show how this spectral
curve arises in an integral representation of the kernel in the case of the matrix model in an external field. We then
show how the expression of this kernel in terms of the spectral curves exhibits a universal behavior in the large matrix
limit depending on the local properties of this curve. We finally apply this procedure to prove the appearance of the Pearcy
and Airy kernels in the context described in the previous sections.

\subsection{From Hermitian matrix integrals to double contour integrals}

In this section, we derive the a double integral representation of the kernel by using the replica formulation introduced by Br\'ezin and
Hikami \cite{Brezin2}. In order to make this paper self-contained, we show
this duality by simply performing gaussian integrals.

Let us consider the partition function
\beq
Z(A):= \int_{H_N} dM e^{-N \Tr \left({M^2\over 2} - A M \right)}
\eeq
where one integrates over hermitian matrices $M$ of size $N \times N$
and $A$ is a deterministic diagonal matrix, with arbitrary $k$\footnote{The previous section considers the particular case $k=2$.}, of the form
\beq
A := \hbox{diag}(\overbrace{a_1, \dots,a_1}^{n_1},\overbrace{a_2,\dots,a_2}^{n_2}, \dots ,\overbrace{a_k,\dots,a_k}^{n_k}).
\eeq

Diagonalizing the matrix $M$ and using the HCIZ formula \cite{HC,IZ}, one is left with the integration over the eigenvalues
$(x_1,x_2,\dots,x_N)$ of $M$:
\beq
Z(A) = \int_{\mathbb{R}^N} \prod_{i=1}^N dx_i {\Delta(x)\over \Delta(a)} e^{-N {\displaystyle \sum_i} \left({x_i^2 \over 2} - x_i a_i\right)}.
\eeq

In order to compute the density of state $R_1(\lambda)$ and the $\ell$-point correlation functions $R_\ell$ defined by
\beq
R_k(\lambda_1,\lambda_2, \dots,\lb_\ell) = {1 \over N^\ell} \left<\displaystyle{\prod_{i=1}^\ell} \Tr \delta(\lambda_i {\mathbb I} -M)\right>,
\eeq
where the average is taken with respect to the probability measure
\beq
{1 \over Z(A)} \prod_{i=1}^N dx_i {\Delta(x)\over \Delta(a)} e^{-N {\displaystyle \sum_i} \left({x_i^2 \over 2} - x_i a_i\right)},
\eeq
we consider their "Fourier" transforms
\beq
U_l(t_1,t_2, \dots,t_l) = \left<\displaystyle{\prod_{i=1}^l} \Tr e^{iNt_iM}\right>
\eeq
and, in particular, one gets the Fourier transform of the two points correlation function:
\beq
U_2(t_1,t_2) = {1 \over Z(A) N^2} \sum_{\alpha_1,\alpha_2=1}^N \int \left(\prod_{j=1}^N dx_j\right) {\Delta(x) \over \Delta(a)}
e^{-N {\displaystyle \sum_{j=1}^N} \left[ {x_j^2 \over 2} - x_j \left(a_j+ i t_1 \delta_{j,\alpha_1}+ i t_2 \delta_{j,\alpha_2}\right)\right]}.
\eeq
One can now integrate the variables $x_j$ by noting that
\beq
\int \left(\prod_{j=1}^N dx_j\right) \Delta(x) e^{-N {\displaystyle \sum_{j=1}^N} \left[ {x_j^2 \over 2} + x_j b_j \right]} = \Delta(b) e^{{N \over 2} {\displaystyle \sum_{j=1}^N} b_j^2}
\eeq
and expanding $\Delta(x) = \prod_{i \neq j} (x_i-x_j)$:
\bea
\lefteqn{U_2(t_1,t_2)}\no\\
& =&\sum_{\alpha_1,\alpha_2=1}^N  e^{N \left(it_1 a_{\alpha_1} + it_2 a_{\alpha_2} - {t_1^2 + t_2^2 \over 2} - t_1 t_2 \delta_{\alpha_1,\alpha_2}\right)} \times \no\\
&& \times { \displaystyle{\prod_{1\leq l < m \leq N}} (a_l-a_m+ i t_1(\delta_{l,\alpha_1}-\delta_{m,\alpha_1})
+ i t_2(\delta_{l,\alpha_2}-\delta_{m,\alpha_2})) \over \displaystyle{\prod_{1\leq l < m \leq N}} (a_l-a_m)}
 .
\eea
One can see that this can be written as a double contour integral
\bea
\lefteqn{U_2(t_1,t_2)}\no\\
& =&{e^{-N {t_1^2 + t_2^2 \over 2}} \over t_1 t_2} \oint \oint {du dv \over (2i \pi )^2} e^{Ni(t_1u+t_2v)}
{(u-v+it_1-i t_2)(u-v) \over (u-v+it_1)(u-v-it_2)} \times \no\\
&& \;\;\;\;\;\;\;\;\;\;\;\;\;\;\;\;\;\;\;\; \times \prod_k \left(1+{it_1 \over u-a_k}\right)\left(1+{it_2 \over v-a_k}\right)\no\\
&=& {e^{-N {t_1^2 + t_2^2 \over 2}} \over t_1 t_2} \oint \oint {du dv \over (2i \pi )^2} e^{Ni(t_1u+t_2v)}
\left(1-{ t_1 t_2 \over (u-v+it_1)(u-v-it_2)}\right) \times \no\\
&& \;\;\;\;\;\;\;\;\;\;\;\;\;\;\;\;\;\;\;\; \times  \prod_k \left(1+{it_1 \over u-a_k}\right)\left(1+{it_2 \over v-a_k}\right),
\eea
where the integration contours encircle all the eigenvalues $a_k$ and the pole $v = u -i t_1$\footnote{see for instance \cite{Brezin2} for more details around eq.(2-20) and eq.(4-40).}.

We can now go back to the correlation function
\beq
R_2(\lambda,\mu) = \int_{-\infty}^\infty \int_{-\infty}^\infty {dt_1 dt_2 \over 4 \pi^2} e^{-iN(t_1\lambda+t_2\mu)} U(t_1,t_2) .
\eeq
By first integrating on $t_1$ and $t_2$ with the shifts $t_1 \to t_1 - iu$ and $t_2 \to t_2 - iu$, this latter equation reads:
\beq
R_2(\lambda,\mu) = K_N(\lambda,\lambda) K_N(\mu,\mu) - K_N(\mu,\lambda) K_N(\lambda,\mu)
\eeq
where the kernel is defined by
\beq
K_N(\lambda,\mu) = \int {dt \over 2 \pi} \oint {dv \over 2i \pi} \displaystyle \prod_{k=1}^N \left({it-a_k \over v-a_k}\right) {1 \over v-it}
e^{-N\left({v^2+t^2 \over 2}+it\lambda-v\mu\right)}
\eeq
where the integration contour for $v$ goes around all the points $a_k$ and the integration for $t$ is parallel to the
real axis and avoids the $v$ contour. Moreover, it can be derived in a very similar way that  any $k$-point function can be written as the Fredholm determinant:
\beq
R_k(x_1,\dots,x_k) = \det \left[K_N(x_i,x_j)\right]_{i,j=1}^k.
\eeq
By Wick rotating the integration variable $t \to i t$, one gets
\beq
K_N(\lambda,\mu) = \int {dt \over 2 i \pi} \oint {dv \over 2i \pi} \displaystyle \prod_{k=1}^N \left({t-a_k \over v-a_k}\right) {1 \over v-t}
e^{-N\left({v^2-t^2 \over 2}+t\lambda-v\mu\right)}
\eeq
where the integration contour for $t$ is now parallel to the imaginary axis.
One can then rewrite it under a more factorized form:
\beq
K_N(\lambda,\mu) = \int {dt \over 2 i \pi} \oint {dv \over 2i \pi} e^{-N(S(\mu,v)-S(\lambda,t))} {1 \over v-t}
\eeq
where
\beq
S(x,y) = {y^2 \over 2} - x y + \sum_{i=1}^k \epsilon_i \ln (y-a_i)
\label{S}\eeq
with the "fraction numbers" given by
\beq
\epsilon_i:=n_i/N.
\label{fraction numbers}\eeq

The first step in the steepest descent analysis of this kernel is to look for the stationary points of the exponent in
the integrand, i.e. we look for  $y$ as a function of $x$ solution of the equation
\beq
\partial_y S(x,y) = y-x + \sum_{i=1}^k {\epsilon_i \over y-a_i} = 0
\label{S'}\eeq
{\em which is nothing but the equation of the classical spectral curve} introduced in the preceding sections (see \ref{F-function}) and in the general study
of the one matrix model in an external field.

\br
If the external matrix $A$ is highly degenerated, the $\epsilon_i$ are fixed and do not depend in $N$. Thus,
the action $S(x,y)$ does not depend on $N$ in this case.
\er

\br
The density of states is given by $\rho(\lambda) = K_N(\lambda,\lambda)$ and its derivative wrt $\lambda$ can be factorized:
\bea
{1 \over N}{\partial \rho(\lambda) \over \partial \lambda} &=&  \int {dt \over 2 i \pi} \oint {dv \over 2i \pi} e^{-N(S_(\lambda,v)-S(\lambda,t))} 
=\phi(\lambda) \psi(\lambda) 
\eea
where
\beq
\phi(\lambda) = \int_{i\mathbb{R}} {i dt \over 2 \pi} e^{N S(\lambda,t)}
~~\mbox{and}~~
\psi(\lambda) = \oint {dv \over 2i \pi} e^{-NS(\lambda,v)}.
\eeq

\er

\subsection{Saddle points, spectral curve and universality}

Let us now forget about the matrix model and consider a general kernel of the form
\beq\label{kernelgeneral}
K_N(\lambda,\mu) = (-1)^N \int {i dt \over 2 \pi} \oint {dv \over 2i \pi} e^{-N(S(\mu,v)-S(\lambda,t))} {1 \over v-t}
\eeq
where the derivative of the action $S(x,y)$:
\beq
\partial_y S(x,y) = {{\cal{E}}(x,y) \over D(x,y)}
\eeq
is a rational function which can be written as the ratio of two polynomials in both variables. Then the locus of stationary points of $S$ is given by an algebraic equation
\beq
{\cal{E}}(x,y)=0,
\eeq
referred to as the {\em spectral curve} in the sequel. Let us study some of its properties necessary to classify the different universal behavior of $K_N(\lambda,\mu)$ as $N \to \infty$.

Generically, the equation
${\cal{E}}(x,y)=0$
has $d_y$ distinct solutions in $y$ for a given value of $x$. Let us denote them $Y_i(x)$ as functions of $x$:
\beq
{\cal{E}}(x,y)=g_{d_y}(x) \prod_{i=1}^{d_y} \left(y-Y_i(x)\right)
\eeq
where $g_{d_y}(x)$ is the leading coefficient of ${\cal{E}}(x,y)$ as a polynomial in $y$.
However, there exists finitely many branch points $x_i$ such that ${\cal{E}}(x,y)=0$ has a double zero, i.e.
two solutions $Y_j(x_i)= Y_l(x_i)$ coincide. One can also characterize them by the property:
\beq
\partial_y\left.{\cal{E}}(x_i,y)\right|_{y:=Y_j(x_i)}=0.
\eeq
More generally, a $\ell$th order branch point
$x_i^{(\ell)}$ is defined by
\beq
 \partial_y^m\left.{\cal{E}}(x_i^{(\ell)},y)\right|_{y:=Y_j(x_i^{(l\ell)})}=0, \mbox{for}~~ m\leq \ell,~~
~~\mbox{with}~~\partial_y^{\ell+1}\left.{\cal{E}}(x_i^{(\ell)},y)\right|_{y:=Y_j(x_i^{(\ell)})}\neq 0.
\eeq
Other types of singularities might occur, but we will not refer to them in this paper.

\br
A $(\ell+1)$th order branch point can be obtained when a $\ell$th order branch point and a simple branch point merge.
\er

In the next sections, we show that, as $N \to \infty$, the kernel $K_N(\lambda,\mu)$ has a universal behavior for
$\lambda$ and $\mu$ approaching the same point with an appropriate scaling (depending on $N$ and the singular behavior
of the spectral curve at this point). We study the first singularities
and show this universal behavior by a {\em local analysis} of the spectral curve whereas the usual Riemann-Hilbert study \cite{Bleher}
involves a global analysis of the latter.

\subsection{Steepest descent analysis}\label{sda}

The asymptotic of this kernel when $N \to \infty$ exhibits different regimes depending on the value of its argument. More
precisely, we get a universal kernel associated to the neighborhood of any point $x_0$ of the spectral curve: if $x_0$
is a generic point of the spectral curve, we get the sine kernel; if $x_0$ is a simple branch point, we get the Airy kernel;
if $x_0$ is a higher order singularity of the curve, we get a universal kernel associated to this particular singularity.

In any case, we proceed with the same steps: we choose a point $x_0$. We then focus on its neighborhood by a change
of variable consistent with the singularity of the curve at $x_0$. We expand the exponent in the kernel around this point
by simply writing a Taylor series expansion.

%
%

\subsubsection{Simple branch point: the Airy kernel}

Consider a one-parameter family of algebraic functions $S(x,y|t)$, parametrized by $t$, such that there exists a critical point $(x_c,y_c,t_c)$ satisfying
\beq
S_y(x_c,y_c|t_c)=S_{yy}(x_c,y_c|t_c)=0,
\eeq
\beq
S_{yyy}(x_c,y_c|t_c)\neq 0
,~~
S_{yx}(x_c,y_c|t_c)\neq 0
,~~
S_{yt}(x_c,y_c|t_c)\neq 0.
\eeq
It means that the spectral curve ${\cal{E}}(x,y|t_c)$ has a simple branch point at $(x_c,y_c)$.
One shows that in the neighborhood of this branch point, the kernel can be rescaled in such a way that it converges to the
Airy kernel as $N \to \infty$.
Consider the changes of variables allowing to focus on the neighborhood of the critical point:
\beq
\left\{
\begin{array}{l}
t:= t_c + {\alpha_t T \over N^{1 \over 3}} \cr
x:= x_c + {\alpha_x  \over N^{1 \over 3}} + {\beta_x X \over N^{2 \over 3}} \cr
y:= y_c + {\alpha_y Y \over N^{1 \over 3}} \cr
\end{array}
\right. ,
\label{7.27}\eeq
and expand $S(x,y|t)$ around this critical point using (\ref{7.27})  as $N \to \infty$.
One then expands $S(x,y|t)$ and $S(\tilde x,\tilde y|t)$ in aTaylor series in $t$, $x$ and $y$, thus obtaining for the expression in the integrand of the kernel:
\bea
\lefteqn{ N \left( S(x,y;t) - S(\tilde{x},\tilde{y};t)\right)}\no\\
&=& N^{1 \over 3} \left(\alpha_y \left[\alpha_x S_{xy}(x_c,y_c;t_c)+ \alpha_t T S_{ty}(x_c,y_c;t_c)\right] (Y - \widetilde{Y})
+ \beta_x X S_{x}(x_c,y_c;t_c) \right)\cr
&& + (Y^3-\widetilde{Y}^3) {\alpha_y^3 \over 6} S_{yyy}(x_c,y_c;t_c) \no\\
&&+ (Y^2-\widetilde{Y}^2) {\alpha_y^2 \over 2} \left[ \alpha_x S_{xyy}(x_c,y_c;t_c)
+ \alpha_t T S_{tyy}(x_c,y_c;t_c) \right] \cr
&& + (XY - \widetilde{X}\widetilde{Y}) \beta_x \alpha_y S_{xy}(x_c,y_c;t_c) \no\\
&&+ (X-\widetilde{X}) \beta_x \left[ \alpha_x  S_{xx}(x_c,y_c;t_c)
+\alpha_t  T S_{tx}(x_c,y_c;t_c) \right]\cr
&& + (Y-\widetilde{Y}) \alpha_y \left[ {\alpha_x^2 \over 2} S_{xxy}(x_c,y_c;t_c) + \alpha_x \alpha_t T S_{txy}(x_c,y_c;t_c)
+ {\alpha_t^2 \over 2} T^2 S_{tty}(x_c,y_c;t_c)\right]\cr
&& + O(N^{-{1 \over 3}});
\eea
which we write for brevity as
\bea
\lefteqn{N \left( S(x,y;t) - S(\tilde{x},\tilde{y};t)\right)=} \no\\
&=& N^{1 \over 3}\left[\alpha_1 \left(Y-\widetilde{Y}\right) + \alpha_2 \left(X-\widetilde{X}\right)\right]
+ \alpha_3 \left(Y^3 - \widetilde{Y}^3 \right) + \alpha_4 \left(Y^2-\widetilde{Y}^2\right)\no\\
&& \quad + \alpha_5 \left(XY - \widetilde{X}\widetilde{Y}\right) + \alpha_6 \left(Y- \widetilde{Y}\right)
+ \alpha_7 \left(X- \widetilde{X}\right)
\eea
for coefficients $\left\{\alpha_i\right\}_{i=1}^7$ functions of the scaling parameters $\alpha_t$, $\alpha_x$, $\beta_x$ and $\alpha_y$.

The coefficients $\alpha_2$ and $\alpha_7$ can be eliminated by conjugation of the kernel leaving the Fredholm
determinants invariant. In order to recover the Airy kernel, one has to fix the remaining coefficients by
\beq
\alpha_1 = 0 \qquad , \qquad  \alpha_3= {1 /3} \qquad , \qquad \alpha_5 = -1
\qquad \hbox{and} \qquad
\alpha_6-\alpha_4^2 = 0.
\eeq
$\alpha_y$ and $\beta_x$ are determined by the constraints on $\alpha_3$ and $\alpha_5$ respectively, while the first
equation gives ${\alpha_t T \over \alpha_x}$ and the last one determines $\alpha_x^2$.

\subsubsection{Double branch point: the Pearcey kernel}\label{secpearcy}

Let us now consider an algebraic a function $S(x,y|t)$ such that there exists a critical point $(x_c,y_c,t_c)$ satisfying
\beq
S_y(x_c,y_c|t_c)=S_{yy}(x_c,y_c|t_c)=S_{yyy}(x_c,y_c|t_c)=0,
\eeq
\beq
S_{yyyy}(x_c,y_c|t_c)\neq 0
,~~
S_{yx}(x_c,y_c|t_c)\neq 0
,~~
S_{yt}(x_c,y_c|t_c)\neq 0.
\eeq
It means that the spectral curve ${\cal{E}}(x,y|t_c)$ has a double branch point at $(x_c,y_c)$.
In the neighborhood of this critical point, one can rescale the kernel in such a way that it converges to the
Pearcey kernel as $N \to \infty$. To this effect,
consider the changes of variables in the neighborhood of the critical point:
\beq
\left\{
\begin{array}{l}
t:= t_c + {\alpha_t T \over N^{1 \over 2}} \cr
x:= x_c + {\alpha_x  \over N^{1 \over 2}} + {\beta_x X \over N^{3 \over 4}} \cr
y:= y_c + {\alpha_y Y \over N^{1 \over 4}} \cr
\end{array}
\right.
\label{7.33}\eeq
and expand $S(x,y|t)$ around the critical point using (\ref{7.33}) as $N \to \infty$.
\bea
&&N \left( S(x,y;t) -  S(\tilde{x},\tilde{y};t)\right) = \cr
&&= N^{1 \over 4} \left[ \beta_x  S_{x}(x_c,y_c;t_c)( X - \widetilde{X}) + \alpha_y (Y- \widetilde{Y}) \left(\alpha_x   S_{xy}(x_c,y_c;t_c) + \alpha_t  T  S_{ty}(x_c,y_c;t_c)\right)\right]\cr
&& + {\alpha_y^4  \over 24} S_{yyyy}(x_c,y_c;t_c) (Y^4 - \widetilde{Y}^4)
\no\\
&&+ {\alpha_y^2 \over 2} \left[\alpha_t  T   S_{tyy}(x_c,y_c;t_c)+ \alpha_x   S_{xyy}(x_c,y_c;t_c)\right] (Y^2 - \widetilde{Y}^2) \cr
&& + \beta_x \alpha_y  S_{xy}(x_c,y_c;t_c) (YX - \widetilde{Y}\widetilde{X}) + O(N^{-{1 \over 4}}).
\eea

\br
One can see that (as in the previous case) all the terms  in this expression, except the scaling factor $\beta_x  S_{x}(x_c,y_c;t_c)( X - \widetilde{X})$,
depend only on derivatives of $S$ wrt $y$. This shows that the kernel depends only on the spectral curve
${\cal{E}}(x,y|t) := S_y (x,y|t)$.
\er

One can see that the variables to be integrated $Y$ and $\widetilde{Y}$ appear only in terms which do not
blow up as $N \to \infty$ except one which is proportional to $N^{1 \over 4}$. One can get rid of this term by fine-tuning
the coefficients of the change of variable. Indeed, imposing the constraint on ${\alpha_t T \over \alpha_x}$,
\beq
\alpha_x   S_{xy}(x_c,y_c;t_c) + \alpha_t  T  S_{ty}(x_c,y_c;t_c) = 0 
\eeq
eliminates this term.
One can finally normalize the scaling coefficients $\alpha_y$, $\alpha_x$ and $\beta_x$ respectively by setting
\beq
{\alpha_y^4 \over 24} S_{yyyy}(x_c,y_c;t_c) = - {1 \over 4},
\qquad \beta_x \alpha_y  S_{xy}(x_c,y_c;t_c)= - 1
\eeq
\beq
{\alpha_x \alpha_y^2 \over 2} \left[{\alpha_t  T \over \alpha_x}   S_{tyy}(x_c,y_c;t_c)+   S_{xyy}(x_c,y_c;t_c)\right] =  {\tau \over 2}
,\eeq
yielding the quartic exponent appearing in the Pearcey kernel:
\bea
\lefteqn{N \left( S(x,y;t) -  S(\tilde{x},\tilde{y};t)\right)}
\no\\
&=& N^{1 \over 4} \beta_x  S_{x}(x_c,y_c;t_c)( X - \widetilde{X}) \no\\
&& + {1 \over 4}(Y^4 - \widetilde{Y}^4)
- {\tau \over 2} (Y^2 - \widetilde{Y}^2)
 + (YX - \widetilde{Y}\widetilde{X}) + O(N^{-{1 \over 4}}),
\eea
where the term in $(X - \widetilde{X})$ can be eliminated by conjugation of the kernel.

\subsubsection{$k$th order branch point}
More generally, consider a point $(x_c,y_c;t_c)$ of the spectral curve ${\cal{E}}(x_c,y_c;t_c)$ where\footnote{Remark that the Airy case corresponds
to $l=1$ and the Pearcey case to $l=2$.}
\beq
 \partial_y^m {\cal{E}}(x_c,y_c;t_c)=0,~~\mbox{for all} ~
m\leq l , ~~\mbox{with}~\partial_y^{l+1} {\cal{E}}(x_c,y_c;t_c)\neq 0,
\eeq
while still assuming that $\partial_t {\cal{E}}$ and $\partial_x {\cal{E}}$ do not vanish at this critical point.

Let us use this example to explain how one can guess the rescaling to obtain the desired universal kernel.
Generically, one looks for a rescaling of the form:
\beq
\left\{
\begin{array}{l}
t:= t_c + {\alpha_t T \over N^{\gamma_t}} \cr
x:= x_c + {\alpha_x X \over N^{\gamma_x}} \cr
y:= y_c + {\alpha_y Y \over N^{\gamma_y}} \cr
\end{array}
\right.
\eeq
for some critical exponents $\gamma_t$, $\gamma_x$ and $\gamma_y$ to be determined. One then writes the Taylor expansion
of the action $S$ around the critical point using this rescaling: this expansion can be seen as a series in $N^{-\gamma}$ for
some exponent $\gamma$ depending on the rescaling. Since the action is multiplied by $N$ in the kernel and 
integrating over the variable $y$, one needs that all the terms depending on $Y$ in this expansion are of order ${1 \over N}$
at most as $N \to \infty$. In particular, considering the Taylor series with respect to $y$:
\beq
S(x_c,y;t_c)= S(x_c,y_c;t_c) +\dots + {\left(\alpha_y Y\right)^{l+2}\over (l+2)!} \partial_y^{l+1} {\cal{E}}(x_c,y_c;t_c) N^{-(l+2)\gamma_y}
+o\left(N^{-(l+2)\gamma_y}\right),
\no\eeq
one must impose
$
(l+2)  \gamma_y = 1
$
for this contribution to appear, thus fixing the critical exponent for $y$:
\beq
\gamma_y = {1 \over l+2}.
\eeq
On the other hand, one expects $x$ and $y$ to couple in the exponent of the limiting integrand (otherwise, it would only give
a simple multiplicative factor). Thus, the first contribution of a mixed derivative wrt to $x$ and $y$ of the action
must be of order ${1 \over N}$ which imposes:
\beq
\gamma_x= 1-\gamma_y = {l+1 \over l+2}.
\eeq

One can remark that we could stop the procedure here and not rescale the time which could remain decoupled from $x$ and $y$.
But one could also try to couple it with $y$ by taking $\alpha_t\neq 0$. If one wants this coupling to be different from the
one between $x$ and $y$, one must impose $\gamma_t<\gamma_x$ and one gets, in the simplest case\footnote{This means that one wants to have a coupling
$T Y^2$, but one could also obtain in a more complicated way higher order times by looking for couplings of the
form $TY^k$ with $k>1$.}
\beq
\gamma_t = 1- 2 \gamma_y = {l \over l+2}.
\eeq
By doing so, one has coupled $t$ to $y$. But one has also introduced a divergent term coming from the coefficient
of $S_{ty}$. One can compensate this term by completing the change of variable for $x$ and writing:
\beq
\left\{
\begin{array}{l}
t:= t_c + {\alpha_t T \over N^{l\over l+2}} \cr
x:= x_c + {\beta_x \over N^{l \over l+2}} + {\alpha_x X \over N^{l+1 \over l+2}} \cr
y:= y_c + {\alpha_y Y \over N^{1 \over l+2}} \cr
\end{array}
\right. .
\eeq

Doing so, one finally gets the associated  kernel in terms of the rescaled variables:
\beq
K_k(X,X') = \int dY \int dY' e^{{Y^l\over l!}- {T Y^2\over 2} + X Y} e^{-{Y'^l\over l!}- {T Y'^2\over 2} + X' Y'} {1 \over Y-Y'}.
\eeq

\br
One omits a detailed proof here, although the details can be filled in. Indeed, in order to do so, one carefully studies the integration contours and the normalization,
which implies a rescaling of the kernel as well as conjugation. 
It is done in the preceding sections for the Airy and Pearcey cases for two ending points.
\er

\subsection{General case}

In a more general case, i.e. the two arguments $x$ and $y$ of $K_N$ approaching any singular or non-singular point of the
algebraic curve, one can use the same study leading to a different universal limiting kernel associated to each type of singularity.
Let us summarize the procedure one has to follow to obtain the universal kernel associated to a given singularity:
\begin{itemize}
\item First compute the multiple derivatives of the action $S(x,y;t)$ with respect to $x$, $y$ and $t$ at the considered
critical point and fix which are the first non-vanishing derivatives;

\item Fix the critical exponent in the rescaling by studying the large $N$ behavior of the Taylor expansion
of the action around the critical point, under consideration;

\item Normalize the change of variables and compute the rescaling of the kernel through the Taylor expansion
of the action of the action around the critical point: one then gets the universal limiting kernel;

\item Finally, carefully study the integration contours to check that they give a right path for the steepest descent analysis.

\end{itemize}

\subsection{Example: back to the matrix model and non-intersecting Brownian motions}

Let us now apply the preceding procedure to the case of the random matrix model in an external field studied in the preceding sections: one
considers $N$ Brownian motions starting from 0 and going to two ending points at $t=1$.

As reminded in section \ref{section1} of this paper, the kernel, given in \eq{BaikKernel},  
can easily be taken to the form of \eq{kernelgeneral} :
\beq
H_n(\bar\mu,\bar\lambda;t)= \int {i du \over 2 \pi} \oint {dv \over 2i \pi} e^{-n(S(\bar\mu,v;t)-S(\bar\lambda,u;t))} {1 \over v-u}
\eeq
by using the rescaling
\beq
v:={V t \over c \sqrt{n}}\, , \; u:={U t \over c \sqrt{n}}\, , \; \tilde{a}:={a t \over c \sqrt{n}}\, , \; \tilde{b}:={b t \over c \sqrt{n}}\, , \;
 \mu:={x \over \sqrt{n}} \; \hbox{and} \; \lambda:={y \over \sqrt{n}}
\eeq
with ${t \over c} = \sqrt{2 t \over 1-t}$ and the action
\beq
S(x,u;t):= {u^2 \over 2} - x {u \over c} + p \log(u-\tilde{a}) + (1-p) \log(u-\tilde{b})
\eeq
which is nothing but the function $S(u,x;t)=F(u)$ of \eq{F-function} studied in section \ref{section2}. The study of this section, and particularly
\eq{derivFnul}, states that
\beq
S_u(x_0,u_0;t_0)=S_{uu}(x_0,u_0;t_0)=S_{uuu}(x_0,u_0;t_0)=0
\eeq
at the critical branch point $(x_0,u_0)$ of the spectral curve
\beq
S_u(x,u;t_0)= u- {x  \over c_0} + {p \over u-\tilde{a}_0 } + {1-p \over u-\tilde{b}_0} =0
\eeq
where one notes $c_0 = c(t_0)$, $\tilde{a}_0 = {a t_0 \over c_0}$ and $\tilde{b}_0 = {b t_0 \over c_0}$.
Moreover, using the notations and computations of section \ref{section2}, one has:
\beq
S_{ux}(x_0,u_0;t_0)= - {1 \over c_0},~~~~{S_{uuuu}(x_0,u_0;t_0)\over 4!} = - {q^2-q+1\over 4q}
\eeq
and
\beq
S_{ut}(x_0,u_0;t_0)= {x_0 (1-2t_0)  \over 2 c_0 t_0 (1-t_0)} + {p \tilde{a}_0 c_0^2 \over t_0^2 (1-t_0)^2\left(u-\tilde{a}_0\right)^2}
+ {(1-p) \tilde{b}_0 c_0^2 \over t_0^2 (1-t_0)^2(u-\tilde{b}_0)^2}
.\eeq
This implies we are in the case of a double branch point studied in section \ref{secpearcy} and thus, with the right rescaling,
the kernel converges to the Pearcey kernel as $n \to \infty$. Let us now check the conditions to obtain precisely the Pearcey
kernel. From the study of section \ref{secpearcy}, one must consider the rescaling
\beq
\left\{
\begin{array}{l}
t:= t_0 + {\alpha_t \tau \over N^{1 \over 2}} \cr
x:= x_0 + {\alpha_x  \over N^{1 \over 2}} + {\beta_x X \over N^{3 \over 4}} \cr
u:= u_0 + {\alpha_u Y \over N^{1 \over 4}} \cr
\end{array}
\right.
\eeq
with the conditions
\beq\label{cond1}
\alpha_x   S_{xu}(x_0,u_0;t_0) + \alpha_t  \tau  S_{tu}(x_0,u_0;t_0) = 0
\eeq
\beq\label{cond2}
{\alpha_u^4 \over 24} S_{uuuu}(x_0,u_0;t_0) = - {1 \over 4},
\eeq
\beq\label{cond3}
{\alpha_u^2 \over 2} \left[\alpha_t  \tau   S_{tuu}(x_0,u_0;t_0)+ \alpha_x   S_{uux}(x_0,u_0;t_0)\right] =  {\tau \over 2}
\eeq
and
\beq\label{cond4}
\beta_x \alpha_u  S_{ux}(x_0,u_0;t_0)=-1.
\eeq

Let us check that these conditions agree with the rescalling \eq{rescaling1}. From conditions (\ref{cond2}) and (\ref{cond4}), one gets ($\mu$ was defined in \eq{constants})
\beq
\alpha_u
= {1 \over \mu},~~~\beta_x= 
c_0 \mu,
\eeq
whereas \eq{cond3} and \eq{cond1} give
\beq
\alpha_t= {(1-t_0)^2 c_0^2 \over 2} {1 \over {p \tilde{a}_0 \over (u_0-\tilde{a}_0)^3} + {(1-p) \tilde{b}_0 \over (u_0-\tilde{b}_0)^3}}
\left({q^2-q+1 \over q}\right)^{1 \over 2}= 2 c_0^2 \mu^2
\eeq
and \eq{cond1} states that
\bean
\alpha_x&=& \tau
\left[{x_0  (1-t_0)(1-2t_0)\over 2 c_0} + {p \tilde{a}_0 c_0^2\over t_0 (u_0-\tilde{a}_0)^2} + {(1-p) \tilde{b}_0 c_0^2\over t_0 (u_0-\tilde{b}_0)^2}\right]
\\
&&~~~~\times {\mu^2 \over {p \tilde{a}_0 \over (u_0-\tilde{a}_0)^3} + {(1-p) \tilde{b}_0 \over (u_0-\tilde{b}_0)^3}}
.
\eean

It is then easily checked by plugging in the values of $x_0$, $t_0$ and $u_0$ in terms of $q$ (with Maple for example) that
this rescaling coincides with \eq{rescaling1}, that is to say, the rescaling considered in the preceding part.

\subsection{Application: Brownian bridges from one point to k points}

We can also use this analysis to study more general statistical systems. Let us now consider $N$ non-intersecting Brownian
motions stating from $0$ at time $t=0$ and ending at $k$ points $a_i$ by groups of $n_i$ particles where
$n_i = \epsilon_i N$. It is a classical result that their probability measure at a given time $t$ is given by the probability
measure of the eigenvalues of an Hermitian random matrix in an external field  with external matrix $A(t)$ whose eigenvalues
are given by $a_i(t) = \sqrt{2 t \over 1-t} a_i$.

Then, the spectral curve is given by:
\beq
{\cal{E}}(x,y) = y-x + \sum_{i=1}^k {\epsilon_i \over y-a_i(t)} = 0
\eeq
 and has genus zero. One can straightforwardly (just solving the equation in $x$) find a rational parametrization under the form:
\beq
x(z) = z + {\displaystyle \sum_{i=1}^k} {\epsilon_i \over  (z-a_i(t))} ,~~~~
y(z) = z 
\eeq
and the branch points are given by
\beq\label{eqbp}
x'(z) = 1 - \sum_{i=1}^k {\epsilon_i \over  (z-a_i(t))^2} = 0.
\eeq

For a generic time $t$, one has $2k$ branch points, $2l$ of them $\left\{z_{2i-1},z_{2i}\right\}_{i=1}^l$ being real and lying in the so-called physical sheet
of the spectral curve (see \cite{AptBleKui,BleKui2,Obrown} for an extensive study of the spectral curve). Thus, by he study of section \ref{section2}
and \ref{sda}, one can conclude that the kernel converges to the Airy kernel in the neighborhood of these real branch points $z_i$ recovering the results of
\cite{Obrown}.

Now, as the time decreases from 1 to 0, some real branch points $z_{2i}(t)$ and $z_{2i+1}(t)$ merge for some critical
value of the time $t_c$. In the Brownian motion setting, it correspond to the times when one big group of particles splits
into two smaller. In terms of the family of spectral curves parameterized by the time $t$, it corresponds to the merging
of two cuts. In the neighborhood of this double branch point, thanks to the study of section \ref{sda}, one can rescale the kernel
in such a way that it converges to the Pearcey kernel as $N \to \infty$. This generalizes the result of th. \ref{Theo 1}
to any cusp in this kind of processes and recovers the former results of \cite{Obrown}.

It is then natural to ask the question:
Can we have higher order singularity in these processes?
When the $\epsilon_i$'s are independent of $N$ and the $a_i$'s real, the answer is no.

Indeed, it amounts to knowing what is the highest order possible for a real root of \eq{eqbp}. This problem can be
rephrased as knowing the higher order possible for a real root of the equation
\beq
T = \sum_{i=1}^k {\epsilon_i \over (z-a_i)^2},
\eeq
with the constraint $\sum_i \epsilon_i = 1$
and the rescaled time evolution $T= {2 t \over 1-t} \in [0,\infty]$.

One can prove that the real roots of this equation are at most double. For this purpose, let us follow the evlution
of the roots as $t$ decreases from $1$ to $0$, i.e. for $T$ going from $\infty$ to 0.

For $T$ large, this equation has obviously $2k$ simple real roots $z_i$ located around their $T \to \infty$ value,
i.e. $\left|z_{2i - 1} - a_i\right| \ll 1$ and $\left|z_{2i } - a_i\right| \ll 1$ with $z_{2i-1} < a_i< z_{2i}$.
Now, for any real solution of this equation, one can compute
\beq
{dz_j \over dT} = - {1 \over 2 \sum_i {\epsilon_i \over (z_j-a_i)^3}}.
\eeq
This derivative does not change sign as long as $z_j$ does not cross any $a_i$. It means that $z_{2i}$ (resp.
$z_{2i+1}$) keeps on going from $a_i$ to $a_{i+1}$ (resp. from $a_{i+1}$ to $a_i$) as $T$ decreases, unless it reaches
another real root. Following this process one sees that the two real roots $z_{2i}$ and $z_{2i+1}$ meet for some critical time $T_c$ giving
birth to a double real root of this equation. In this first part of the process, one can thus only encounter double
real roots.

Let us keep on decreasing time. For $T<T_c$ and close to it, the double root gives rise to two simple complex
conjugated roots. Let us thus now consider a simple complex root $z = r + i \theta$. For $\theta$ close to $0$, one
can compute
\beq
{d\theta \over d T} = -{1 \over 6 \sum_i {\epsilon_i \theta \over (r-a_i)^4}}
\eeq
to first order in $\theta$. The complex roots are thus repelled by the real axis and cannot thus give birth to real
roots. In this second part of the process, one do not have real roots anymore. The only multiple real roots are thus
obtained when $a_i<z_{2i} \to z_{2i+1}<a_{i+1}$.

\bibliographystyle{amsalpha}

\end{document}